\documentclass[10pt]{amsart}
\usepackage[english,russian]{babel}
\usepackage[cp1251]{inputenc}
\usepackage{amsmath}
\usepackage{amssymb}
\usepackage{amsfonts}

\usepackage[linktocpage=true, colorlinks=true, linkcolor=blue, citecolor=blue, urlcolor=blue]{hyperref}

\begin{document}
\renewcommand{\refname}{References}
\renewcommand\contentsname{Contents}

\thispagestyle{empty}

\title[Expansion of Iterated Stratonovich Stochastic Integrals]
{Expansion of Iterated Stratonovich Stochastic Integrals of Multiplicity 2
Based on Double Fourier--Legendre Series Summarized by Pringsheim Method}
\author[D.F. Kuznetsov]{Dmitriy F. Kuznetsov}
\address{Dmitriy Feliksovich Kuznetsov
\newline\hphantom{iii} Peter the Great Saint-Petersburg Polytechnic University,
\newline\hphantom{iii} Polytechnicheskaya ul., 29,
\newline\hphantom{iii} 195251, Saint-Petersburg, Russia}%
\email{sde\_kuznetsov@inbox.ru}

\maketitle {\small
\begin{quote}
\vspace{5mm}
\noindent{\sc Abstract.} 
The article is devoted to
the expansion of iterated Stratonovich stochastic integrals of second
multiplicity into the double series of products of 
standard Gaussian random variables.
The proof of expansion is based on the application of double 
Fourier--Legendre series and 
double trigonometric Fourier series summarized by Pringsheim method. 
The results of the article are generalized to the case
of an arbitrary 
complete orthonormal system of functions 
in the space $L_2([t, T])$ and 
$\psi_1(\tau), \psi_2(\tau)\in L_2([t, T]),$
where $\psi_1(\tau), \psi_2(\tau)$ are
weight functions of the iterated Stratonovich
stochastic integral of second multiplicity.
The considered expansion
can be applied to the numerical integration of Ito stochastic differential 
equations. 
Some recent results on the expansion of iterated Stratonovich
stochastic integrals of multiplicities 3 to 6 are given.

\medskip
\end{quote}
}

\vspace{5mm}

\linespread{1.3}

\tableofcontents

\linespread{1.0}

\section{Introduction}

\vspace{5mm}

Let $(\Omega,$ ${\rm F},$ ${\sf P})$ be a complete probability space, let 
$\{{\rm F}_t, t\in[0,T]\}$ be a nondecreasing right-continous 
family of $\sigma$-algebras of ${\rm F},$
and let ${\bf f}_t$ be a standard $m$-dimensional Wiener stochastic 
process, which is
${\rm F}_t$-measurable for any $t\in[0, T].$ We assume that the components
${\bf f}_{t}^{(i)}$ $(i=1,\ldots,m)$ of this process 
are independent. Consider
an Ito stochastic differential equation (SDE) 
in the integral form

\vspace{1mm}
\begin{equation}
\label{1.5.2}
{\bf x}_t={\bf x}_0+\int\limits_0^t {\bf a}({\bf x}_{\tau},\tau)d\tau+
\int\limits_0^t B({\bf x}_{\tau},\tau)d{\bf f}_{\tau},\ \ \
{\bf x}_0={\bf x}(0,\omega).
\end{equation}

\vspace{4mm}
\noindent
Here ${\bf x}_t$ is some $n$-dimensional stochastic process 
satisfying the equation (\ref{1.5.2}). 
The nonrandom functions ${\bf a}: \mathbb{R}^n\times[0, T]\to\mathbb{R}^n$,
$B: \mathbb{R}^n\times[0, T]\to\mathbb{R}^{n\times m}$
guarantee the existence and uniqueness up to stochastic equivalence 
of a solution
of the Ito SDE (\ref{1.5.2}) \cite{1}. The second integral on 
the right-hand side of (\ref{1.5.2}) is 
interpreted as an Ito stochastic integral.
Let ${\bf x}_0$ be an $n$-dimensional random variable, which is 
${\rm F}_0$-measurable and 
${\sf M}\bigl\{\left|{\bf x}_0\right|^2\bigr\}<\infty$ 
(${\sf M}$ denotes a mathematical expectation).
We assume that
${\bf x}_0$ and ${\bf f}_t-{\bf f}_0$ are independent when $t>0.$

It is well known that one of the effective approaches 
to the numerical integration of 
Ito SDEs is an approach based on the Taylor--Ito and 
Taylor--Stratonovich expansions
\cite{Mi2}-\cite{miltret}. Moreover, one of the most important 
features of such 
expansions is a presence in them of the so-called iterated
Ito and Stratonovich stochastic integrals, which play the key 
role for solving the 
problem of numerical integration of Ito SDEs and have the 
following form

\vspace{1mm}
\begin{equation}
\label{sodom20}
J[\psi^{(k)}]_{T,t}=\int\limits_t^T\psi_k(t_k) \ldots \int\limits_t^{t_{2}}
\psi_1(t_1) d{\bf w}_{t_1}^{(i_1)}\ldots
d{\bf w}_{t_k}^{(i_k)},
\end{equation}

\vspace{2mm}

\begin{equation}
\label{605}
J^{*}[\psi^{(k)}]_{T,t}=
{\int\limits_t^{*}}^{T}
\psi_k(t_k) \ldots 
{\int\limits_t^{*}}^{t_2}
\psi_1(t_1) d{\bf w}_{t_1}^{(i_1)}\ldots
d{\bf w}_{t_k}^{(i_k)},
\end{equation}

\vspace{4mm}
\noindent
where every $\psi_l(\tau)$ $(l=1,\ldots,k)$ is a 
nonrandom function 
at the interval $[t,T],$ ${\bf w}_{\tau}^{(i)}={\bf f}_{\tau}^{(i)}$
for $i=1,\ldots,m$ and
${\bf w}_{\tau}^{(0)}=\tau,$ $i_1,\ldots,i_k = 0, 1,\ldots,m,$

\vspace{-1mm}
$$
\int\limits\ \hbox{and}\ \int\limits^{*}
$$ 

\vspace{4mm}
\noindent
denote Ito and 
Stratonovich stochastic integrals,
respectively. In this paper we use the definition of the Stratonovich 
stochastic integral from \cite{KlPl2}.

Note that usually in applications 
the functions $\psi_l(\tau)\equiv 1$ $(l=1,\ldots,k)$
are equal to $1$ or have a binomial form.
More precisely, 
$\psi_l(\tau)\equiv 1$ $(l=1,\ldots,k)$ and
$i_1,\ldots,i_k = 0, 1,\ldots,m$ in
\cite{Mi2}-\cite{KlPl1}. At the same time
$\psi_l(\tau)\equiv (t-\tau)^{q_l}$ ($l=1,\ldots,k$,\ \ 
$q_1,\ldots,q_k=0, 1, 2,\ldots $) and $i_1,\ldots,i_k = 1,\ldots,m$ in
\cite{3}-\cite{32}.

Effective solution 
of the problem of
combined mean-square approximation for collections 
of iterated Stratonovich stochastic integrals 
(\ref{605}) of second multiplicity
composes the subject of this article.

It is well known that 
the mean-square approximation of iterated Ito and Stratonovich 
stochastic
integrals (\ref{sodom20}), (\ref{605}) using 
multiple integral sums
requires significant
computational costs \cite{7} since this approach
implies the partitioning of the integration interval $[t, T]$ of 
the iterated stochastic
integrals (\ref{sodom20}), (\ref{605})
into parts ($T-t$ is already a sufficiently small value since
$T-t$  plays the role of an integration step in numerical methods for 
solving Ito SDEs).

More efficient approximation methods for
the iterated stochastic
integrals (\ref{sodom20}), (\ref{605})
use Fourier series, and 
they do not require the interval $[t, T]$
to be subdivided into smaller parts. One such method was 
proposed in \cite{Mi2} and elaborated in \cite{KlPl2}, \cite{KPW}. 
This method, which received widespread use, is based on 
the Karhunen--Loeve expansion of the
Brownian bridge process \cite{Mi2} in the eigenfunctions of its 
covariance, which form a complete orthonormal
trigonometric basis in the space $L_2([t, T]).$ 

Note that in \cite{7} (2006) 
the more general and effective 
method (the so-called method of 
generalized multiple Fourier series) for
the mean-square approximation of iterated Ito stochastic
integrals (\ref{sodom20}) was proposed.  
This method is based on the generalized
multiple Fourier series that converge in the sense of norm 
in Hilbert space $L_2([t, T]^k)$, 
where $[t, T]^k$ is the hypercube $[t, T]\times \ldots \times 
[t, T]$ ($k$ times) and $k$ is
the multiplicity of the iterated Ito stochastic integral. 
The method of generalized multiple 
Fourier series was developed in \cite{8}-\cite{31aa},
\cite{arxiv-14}-\cite{new-2023a}.

An important feature of the  method of 
generalized multiple Fourier series is that various complete 
orthonormal systems of functions in the space $L_2(t, T])$ can
be used (the method proposed in \cite{Mi2} 
admits only the trigonometric system of functions).
Hence, we can state the problem of comparing the efficiency 
of using different complete orthonormal
systems of functions in the space $L_2(t, T])$ in the context of
numerical solution of Ito SDEs.
This problem has been solved in \cite{29}, \cite{30} (also see 
\cite{20xxxxx}-\cite{12aa-afterxxx}).
In particular, in \cite{20xxxxx}-\cite{12aa-afterxxx}, \cite{29}, \cite{30} 
it was shown that the optimal system of basis functions 
in the framework of 
numerical solution of Ito SDEs
is the system of Legendre 
polynomials. This fact is true at least for high-order strong
numerical methods with orders of convergence $1.5, 2.0,\ldots $
That is why the part of this article is devoted to the expansions
of iterated Stratonovich stochastic integrals with multiplicity 2 based
on multiple Fourier--Legendre series.

Usage of Fourier series with respect to the system of Legendre polynomials 
for approximation of iterated stochastic integrals took place for 
the first time in \cite{3} (1997), \cite{4} (1998), \cite{5} (also see
\cite{6}-\cite{31aa}, \cite{arxiv-14}-\cite{trace}).

The results of \cite{3}
(also see \cite{4}-\cite{31aa},
\cite{arxiv-14}-\cite{100-000-6}) convincingly testify
that there is a doubtless relation between 
multiplier factor $1/2$, 
which is typical for Stratonovich stochastic integral and included 
into the sum connecting Stratonovich and Ito stochastic integrals, 
and the fact that in the point of finite discontinuity of piecewise
smooth function $f(x)$ its Fourier series  
converges to the value 

\vspace{1mm}
$$
\frac{f(x+0)+f(x-0)}{2}.
$$ 

\vspace{3mm}
\noindent
In addition, in \cite{3}, \cite{4},
\cite{12}-\cite{16}, \cite{19}-\cite{21},
\cite{23}, \cite{25}, \cite{27}, \cite{30a},
\cite{arxiv-7}-\cite{arxiv-8ee}, \cite{100-000-1},
\cite{100-000-5} 
several theorems on expansion of iterated 
Stratonovich stochastic integrals 
were formulated and proved.
As shown in these papers, the final formulas for 
expansions of iterated Stratonovich stochastic integrals 
are more compact than their analogues for iterated Ito stochastic 
integrals. 

This paper continues the study of the relationships between
generalized multiple Fourier series and iterated stochastic integrals.
We use the double Fourier--Legendre series
and double trigonometric Fourier series
(summarized by Pringsheim method)
for the proof of Theorem 5.3 \cite{20}
or Theorem 2.1 \cite{20xxxxx} (also see \cite{20eee}, \cite{12aa-afterxxx}). 
As shown below 
the conditions of these theorems can be weakened.

Moreover, the mentioned theorems are generalized to the case of an arbitrary 
complete orthonormal system of functions
in the space $L_2([t, T])$ 
and $\psi_1(\tau), \psi_2(\tau)\in L_2([t, T]).$

\vspace{5mm}

\section{Method of Generalized Multiple Fourier Series}

\vspace{5mm}

Let us consider an approach to the expansion of iterated
Ito stochastic integrals \cite{7}-\cite{31aa},
\cite{arxiv-14}-\cite{new-2023a}
(the so-called
method of generalized
multiple Fourier series).

Suppose that $\psi_1(\tau),\ldots,\psi_k(\tau)\in L_2([t, T]).$ 
Define the following function on the hypercube $[t, T]^k$

\begin{equation}
\label{ppp}
K(t_1,\ldots,t_k)=
\begin{cases}
\psi_1(t_1)\ldots \psi_k(t_k)\ \ &\hbox{for}\ \ t_1<\ldots<t_k\\
~\\
~\\
0\ \ &\hbox{otherwise}
\end{cases},\ \ \ \ \ t_1,\ldots,t_k\in[t, T],\ \ k\ge 2,
\end{equation}

\vspace{4mm}
\noindent
and 
$K(t_1)\equiv\psi_1(t_1),\ t_1\in[t, T].$

Suppose that $\{\phi_j(x)\}_{j=0}^{\infty}$
is a complete orthonormal system of functions in the space
$L_2([t, T])$. 
The function $K(t_1,\ldots,t_k)$ belongs to the space
$L_2([t, T]^k).$
At this situation it is well-known that the generalized 
multiple Fourier series 
of $K(t_1,\ldots,t_k)\in L_2([t, T]^k)$ is converging 
to $K(t_1,\ldots,t_k)$ in the hypercube $[t, T]^k$ in 
the mean-square sense, i.e.

$$
\hbox{\vtop{\offinterlineskip\halign{
\hfil#\hfil\cr
{\rm lim}\cr
$\stackrel{}{{}_{p_1,\ldots,p_k\to \infty}}$\cr
}} }\Biggl\Vert
K(t_1,\ldots,t_k)-
\sum_{j_1=0}^{p_1}\ldots \sum_{j_k=0}^{p_k}
C_{j_k\ldots j_1}\prod_{l=1}^{k} \phi_{j_l}(t_l)\Biggr\Vert_{L_2([t,T]^k)}=0,
$$

\vspace{3mm}
\noindent
where
\begin{equation}
\label{ppppa}
C_{j_k\ldots j_1}=\int\limits_{[t,T]^k}
K(t_1,\ldots,t_k)\prod_{l=1}^{k}\phi_{j_l}(t_l)dt_1\ldots dt_k,
\end{equation}

\vspace{2mm}
$$
\left\Vert f\right\Vert_{L_2([t,T]^k)}=\left(\int\limits_{[t,T]^k}
f^2(t_1,\ldots,t_k)dt_1\ldots dt_k\right)^{1/2}.
$$

\vspace{5mm}

Consider the partition $\{\tau_j\}_{j=0}^N$ of the interval
$[t,T]$ such that

\begin{equation}
\label{1111}
t=\tau_0<\ldots <\tau_N=T,\ \ \
\Delta_N=
\hbox{\vtop{\offinterlineskip\halign{
\hfil#\hfil\cr
{\rm max}\cr
$\stackrel{}{{}_{0\le j\le N-1}}$\cr
}} }\Delta\tau_j\to 0\ \ \hbox{if}\ \ N\to \infty,\ \ \
\Delta\tau_j=\tau_{j+1}-\tau_j.
\end{equation}

\vspace{3mm}

{\bf Theorem 1} \cite{7} (2006) (also see \cite{8}-\cite{31aa},
\cite{arxiv-14}-\cite{new-2023a}). 
{\it Suppose that
every $\psi_l(\tau)$ $(l=1,\ldots, k)$ is a continuous nonrandom function on 
$[t, T]$ and
$\{\phi_j(x)\}_{j=0}^{\infty}$ is a complete orthonormal system  
of continuous functions in $L_2([t,T]).$ Then

\vspace{2mm}
$$
J[\psi^{(k)}]_{T,t}\  =\ 
\hbox{\vtop{\offinterlineskip\halign{
\hfil#\hfil\cr
{\rm l.i.m.}\cr
$\stackrel{}{{}_{p_1,\ldots,p_k\to \infty}}$\cr
}} }\sum_{j_1=0}^{p_1}\ldots\sum_{j_k=0}^{p_k}
C_{j_k\ldots j_1}\Biggl(
\prod_{l=1}^k\zeta_{j_l}^{(i_l)}\ -
\Biggr.
$$

\vspace{3mm}
\begin{equation}
\label{tyyy}
-\ \Biggl.
\hbox{\vtop{\offinterlineskip\halign{
\hfil#\hfil\cr
{\rm l.i.m.}\cr
$\stackrel{}{{}_{N\to \infty}}$\cr
}} }\sum_{(l_1,\ldots,l_k)\in {\rm G}_k}
\phi_{j_{1}}(\tau_{l_1})
\Delta{\bf w}_{\tau_{l_1}}^{(i_1)}\ldots
\phi_{j_{k}}(\tau_{l_k})
\Delta{\bf w}_{\tau_{l_k}}^{(i_k)}\Biggr),
\end{equation}

\vspace{6mm}
\noindent
where $J[\psi^{(k)}]_{T,t}$ is defined by {\rm (\ref{sodom20}),}

$$
{\rm G}_k={\rm H}_k\backslash{\rm L}_k,\ \ \
{\rm H}_k=\bigl\{(l_1,\ldots,l_k):\ l_1,\ldots,l_k=0,\ 1,\ldots,N-1\bigr\},
$$

\vspace{-1mm}
$$
{\rm L}_k=\bigl\{(l_1,\ldots,l_k):\ l_1,\ldots,l_k=0,\ 1,\ldots,N-1;\
l_g\ne l_r\ (g\ne r);\ g, r=1,\ldots,k\bigr\},
$$

\vspace{4mm}
\noindent
${\rm l.i.m.}$ is a limit in the mean-square sense$,$
$i_1,\ldots,i_k=0,1,\ldots,m,$

\vspace{-1mm}
\begin{equation}
\label{rr23}
\zeta_{j}^{(i)}=
\int\limits_t^T \phi_{j}(s) d{\bf w}_s^{(i)}
\end{equation} 

\vspace{3mm}
\noindent
are independent standard Gaussian random variables
for various
$i$ or $j$ {\rm(}if $i\ne 0${\rm),}
$C_{j_k\ldots j_1}$ is the Fourier coefficient {\rm(\ref{ppppa}),}
$\Delta{\bf w}_{\tau_{j}}^{(i)}=
{\bf w}_{\tau_{j+1}}^{(i)}-{\bf w}_{\tau_{j}}^{(i)}$
$(i=0, 1,\ldots,m),$
$\left\{\tau_{j}\right\}_{j=0}^{N}$ is a partition of
the interval $[t, T],$ which satisfies the condition {\rm (\ref{1111})}.
}

\vspace{2mm}

Note that the condition of continuity of the functions $\phi_j(x)$
$(j=0, 1,\ldots )$ can be 
weakened (see \cite{7}-\cite{16}, \cite{19}-\cite{12aa-afterxxx}).
Another versions and generalizations 
of Theorem 1 can be found in \cite{8}-\cite{31aa},
\cite{arxiv-14}-\cite{new-2023a}.

In order to evaluate the significance of Theorem 1 for practice we will
demonstrate its transformed particular cases for 
$k=1,\ldots,5$ \cite{7}-\cite{31aa},
\cite{arxiv-14}-\cite{100-000-6} 

\vspace{2mm}
\begin{equation}
\label{a1}
J[\psi^{(1)}]_{T,t}
=\hbox{\vtop{\offinterlineskip\halign{
\hfil#\hfil\cr
{\rm l.i.m.}\cr
$\stackrel{}{{}_{p_1\to \infty}}$\cr
}} }\sum_{j_1=0}^{p_1}
C_{j_1}\zeta_{j_1}^{(i_1)},
\end{equation}

\vspace{5mm}
\begin{equation}
\label{leto5001}
J[\psi^{(2)}]_{T,t}
=\hbox{\vtop{\offinterlineskip\halign{
\hfil#\hfil\cr
{\rm l.i.m.}\cr
$\stackrel{}{{}_{p_1,p_2\to \infty}}$\cr
}} }\sum_{j_1=0}^{p_1}\sum_{j_2=0}^{p_2}
C_{j_2j_1}\Biggl(\zeta_{j_1}^{(i_1)}\zeta_{j_2}^{(i_2)}
-{\bf 1}_{\{i_1=i_2\ne 0\}}
{\bf 1}_{\{j_1=j_2\}}\Biggr),
\end{equation}

\vspace{8mm}
$$
J[\psi^{(3)}]_{T,t}=
\hbox{\vtop{\offinterlineskip\halign{
\hfil#\hfil\cr
{\rm l.i.m.}\cr
$\stackrel{}{{}_{p_1,\ldots,p_3\to \infty}}$\cr
}} }\sum_{j_1=0}^{p_1}\sum_{j_2=0}^{p_2}\sum_{j_3=0}^{p_3}
C_{j_3j_2j_1}\Biggl(
\zeta_{j_1}^{(i_1)}\zeta_{j_2}^{(i_2)}\zeta_{j_3}^{(i_3)}
-\Biggr.
$$

\vspace{1mm}
\begin{equation}
\label{a3}
\Biggl.-{\bf 1}_{\{i_1=i_2\ne 0\}}
{\bf 1}_{\{j_1=j_2\}}
\zeta_{j_3}^{(i_3)}
-{\bf 1}_{\{i_2=i_3\ne 0\}}
{\bf 1}_{\{j_2=j_3\}}
\zeta_{j_1}^{(i_1)}-
{\bf 1}_{\{i_1=i_3\ne 0\}}
{\bf 1}_{\{j_1=j_3\}}
\zeta_{j_2}^{(i_2)}\Biggr),
\end{equation}

\vspace{7mm}
$$
J[\psi^{(4)}]_{T,t}
=
\hbox{\vtop{\offinterlineskip\halign{
\hfil#\hfil\cr
{\rm l.i.m.}\cr
$\stackrel{}{{}_{p_1,\ldots,p_4\to \infty}}$\cr
}} }\sum_{j_1=0}^{p_1}\ldots\sum_{j_4=0}^{p_4}
C_{j_4\ldots j_1}\Biggl(
\prod_{l=1}^4\zeta_{j_l}^{(i_l)}
\Biggr.
-
$$
$$
-
{\bf 1}_{\{i_1=i_2\ne 0\}}
{\bf 1}_{\{j_1=j_2\}}
\zeta_{j_3}^{(i_3)}
\zeta_{j_4}^{(i_4)}
-
{\bf 1}_{\{i_1=i_3\ne 0\}}
{\bf 1}_{\{j_1=j_3\}}
\zeta_{j_2}^{(i_2)}
\zeta_{j_4}^{(i_4)}-
$$
$$
-
{\bf 1}_{\{i_1=i_4\ne 0\}}
{\bf 1}_{\{j_1=j_4\}}
\zeta_{j_2}^{(i_2)}
\zeta_{j_3}^{(i_3)}
-
{\bf 1}_{\{i_2=i_3\ne 0\}}
{\bf 1}_{\{j_2=j_3\}}
\zeta_{j_1}^{(i_1)}
\zeta_{j_4}^{(i_4)}-
$$
$$
-
{\bf 1}_{\{i_2=i_4\ne 0\}}
{\bf 1}_{\{j_2=j_4\}}
\zeta_{j_1}^{(i_1)}
\zeta_{j_3}^{(i_3)}
-
{\bf 1}_{\{i_3=i_4\ne 0\}}
{\bf 1}_{\{j_3=j_4\}}
\zeta_{j_1}^{(i_1)}
\zeta_{j_2}^{(i_2)}+
$$
$$
+
{\bf 1}_{\{i_1=i_2\ne 0\}}
{\bf 1}_{\{j_1=j_2\}}
{\bf 1}_{\{i_3=i_4\ne 0\}}
{\bf 1}_{\{j_3=j_4\}}
+
{\bf 1}_{\{i_1=i_3\ne 0\}}
{\bf 1}_{\{j_1=j_3\}}
{\bf 1}_{\{i_2=i_4\ne 0\}}
{\bf 1}_{\{j_2=j_4\}}+
$$
\begin{equation}
\label{a4}
+\Biggl.
{\bf 1}_{\{i_1=i_4\ne 0\}}
{\bf 1}_{\{j_1=j_4\}}
{\bf 1}_{\{i_2=i_3\ne 0\}}
{\bf 1}_{\{j_2=j_3\}}\Biggr),
\end{equation}

\vspace{7mm}
$$
J[\psi^{(5)}]_{T,t}
=\hbox{\vtop{\offinterlineskip\halign{
\hfil#\hfil\cr
{\rm l.i.m.}\cr
$\stackrel{}{{}_{p_1,\ldots,p_5\to \infty}}$\cr
}} }\sum_{j_1=0}^{p_1}\ldots\sum_{j_5=0}^{p_5}
C_{j_5\ldots j_1}\Biggl(
\prod_{l=1}^5\zeta_{j_l}^{(i_l)}
-\Biggr.
$$
$$
-
{\bf 1}_{\{i_1=i_2\ne 0\}}
{\bf 1}_{\{j_1=j_2\}}
\zeta_{j_3}^{(i_3)}
\zeta_{j_4}^{(i_4)}
\zeta_{j_5}^{(i_5)}-
{\bf 1}_{\{i_1=i_3\ne 0\}}
{\bf 1}_{\{j_1=j_3\}}
\zeta_{j_2}^{(i_2)}
\zeta_{j_4}^{(i_4)}
\zeta_{j_5}^{(i_5)}-
$$
$$
-
{\bf 1}_{\{i_1=i_4\ne 0\}}
{\bf 1}_{\{j_1=j_4\}}
\zeta_{j_2}^{(i_2)}
\zeta_{j_3}^{(i_3)}
\zeta_{j_5}^{(i_5)}-
{\bf 1}_{\{i_1=i_5\ne 0\}}
{\bf 1}_{\{j_1=j_5\}}
\zeta_{j_2}^{(i_2)}
\zeta_{j_3}^{(i_3)}
\zeta_{j_4}^{(i_4)}-
$$
$$
-
{\bf 1}_{\{i_2=i_3\ne 0\}}
{\bf 1}_{\{j_2=j_3\}}
\zeta_{j_1}^{(i_1)}
\zeta_{j_4}^{(i_4)}
\zeta_{j_5}^{(i_5)}-
{\bf 1}_{\{i_2=i_4\ne 0\}}
{\bf 1}_{\{j_2=j_4\}}
\zeta_{j_1}^{(i_1)}
\zeta_{j_3}^{(i_3)}
\zeta_{j_5}^{(i_5)}-
$$
$$
-
{\bf 1}_{\{i_2=i_5\ne 0\}}
{\bf 1}_{\{j_2=j_5\}}
\zeta_{j_1}^{(i_1)}
\zeta_{j_3}^{(i_3)}
\zeta_{j_4}^{(i_4)}
-{\bf 1}_{\{i_3=i_4\ne 0\}}
{\bf 1}_{\{j_3=j_4\}}
\zeta_{j_1}^{(i_1)}
\zeta_{j_2}^{(i_2)}
\zeta_{j_5}^{(i_5)}-
$$
$$
-
{\bf 1}_{\{i_3=i_5\ne 0\}}
{\bf 1}_{\{j_3=j_5\}}
\zeta_{j_1}^{(i_1)}
\zeta_{j_2}^{(i_2)}
\zeta_{j_4}^{(i_4)}
-{\bf 1}_{\{i_4=i_5\ne 0\}}
{\bf 1}_{\{j_4=j_5\}}
\zeta_{j_1}^{(i_1)}
\zeta_{j_2}^{(i_2)}
\zeta_{j_3}^{(i_3)}+
$$
$$
+
{\bf 1}_{\{i_1=i_2\ne 0\}}
{\bf 1}_{\{j_1=j_2\}}
{\bf 1}_{\{i_3=i_4\ne 0\}}
{\bf 1}_{\{j_3=j_4\}}\zeta_{j_5}^{(i_5)}+
{\bf 1}_{\{i_1=i_2\ne 0\}}
{\bf 1}_{\{j_1=j_2\}}
{\bf 1}_{\{i_3=i_5\ne 0\}}
{\bf 1}_{\{j_3=j_5\}}\zeta_{j_4}^{(i_4)}+
$$
$$
+
{\bf 1}_{\{i_1=i_2\ne 0\}}
{\bf 1}_{\{j_1=j_2\}}
{\bf 1}_{\{i_4=i_5\ne 0\}}
{\bf 1}_{\{j_4=j_5\}}\zeta_{j_3}^{(i_3)}+
{\bf 1}_{\{i_1=i_3\ne 0\}}
{\bf 1}_{\{j_1=j_3\}}
{\bf 1}_{\{i_2=i_4\ne 0\}}
{\bf 1}_{\{j_2=j_4\}}\zeta_{j_5}^{(i_5)}+
$$
$$
+
{\bf 1}_{\{i_1=i_3\ne 0\}}
{\bf 1}_{\{j_1=j_3\}}
{\bf 1}_{\{i_2=i_5\ne 0\}}
{\bf 1}_{\{j_2=j_5\}}\zeta_{j_4}^{(i_4)}+
{\bf 1}_{\{i_1=i_3\ne 0\}}
{\bf 1}_{\{j_1=j_3\}}
{\bf 1}_{\{i_4=i_5\ne 0\}}
{\bf 1}_{\{j_4=j_5\}}\zeta_{j_2}^{(i_2)}+
$$
$$
+
{\bf 1}_{\{i_1=i_4\ne 0\}}
{\bf 1}_{\{j_1=j_4\}}
{\bf 1}_{\{i_2=i_3\ne 0\}}
{\bf 1}_{\{j_2=j_3\}}\zeta_{j_5}^{(i_5)}+
{\bf 1}_{\{i_1=i_4\ne 0\}}
{\bf 1}_{\{j_1=j_4\}}
{\bf 1}_{\{i_2=i_5\ne 0\}}
{\bf 1}_{\{j_2=j_5\}}\zeta_{j_3}^{(i_3)}+
$$
$$
+
{\bf 1}_{\{i_1=i_4\ne 0\}}
{\bf 1}_{\{j_1=j_4\}}
{\bf 1}_{\{i_3=i_5\ne 0\}}
{\bf 1}_{\{j_3=j_5\}}\zeta_{j_2}^{(i_2)}+
{\bf 1}_{\{i_1=i_5\ne 0\}}
{\bf 1}_{\{j_1=j_5\}}
{\bf 1}_{\{i_2=i_3\ne 0\}}
{\bf 1}_{\{j_2=j_3\}}\zeta_{j_4}^{(i_4)}+
$$
$$
+
{\bf 1}_{\{i_1=i_5\ne 0\}}
{\bf 1}_{\{j_1=j_5\}}
{\bf 1}_{\{i_2=i_4\ne 0\}}
{\bf 1}_{\{j_2=j_4\}}\zeta_{j_3}^{(i_3)}+
{\bf 1}_{\{i_1=i_5\ne 0\}}
{\bf 1}_{\{j_1=j_5\}}
{\bf 1}_{\{i_3=i_4\ne 0\}}
{\bf 1}_{\{j_3=j_4\}}\zeta_{j_2}^{(i_2)}+
$$
$$
+
{\bf 1}_{\{i_2=i_3\ne 0\}}
{\bf 1}_{\{j_2=j_3\}}
{\bf 1}_{\{i_4=i_5\ne 0\}}
{\bf 1}_{\{j_4=j_5\}}\zeta_{j_1}^{(i_1)}+
{\bf 1}_{\{i_2=i_4\ne 0\}}
{\bf 1}_{\{j_2=j_4\}}
{\bf 1}_{\{i_3=i_5\ne 0\}}
{\bf 1}_{\{j_3=j_5\}}\zeta_{j_1}^{(i_1)}+
$$
\begin{equation}
\label{a5}
+\Biggl.
{\bf 1}_{\{i_2=i_5\ne 0\}}
{\bf 1}_{\{j_2=j_5\}}
{\bf 1}_{\{i_3=i_4\ne 0\}}
{\bf 1}_{\{j_3=j_4\}}\zeta_{j_1}^{(i_1)}\Biggr),
\end{equation}

\vspace{8mm}
\noindent
where ${\bf 1}_A$ is the indicator of the set $A$.

For further consideration, let us 
consider the generalization of formulas (\ref{a1})--(\ref{a5})                 
for the case of an arbitrary multiplicity $k$ $(k\in\mathbb{N})$ of 
the iterated Ito stochastic integral $J[\psi^{(k)}]_{T,t}$ defined by (\ref{sodom20}).
In order to do this, let us
introduce some notations. 
Consider the unordered
set $\{1, 2, \ldots, k\}$ 
and separate it into two parts:
the first part consists of $r$ unordered 
pairs (sequence order of these pairs is also unimportant) and the 
second one consists of the 
remaining $k-2r$ numbers.
So, we have

\begin{equation}
\label{leto5007}
(\{
\underbrace{\{g_1, g_2\}, \ldots, 
\{g_{2r-1}, g_{2r}\}}_{\small{\hbox{part 1}}}
\},
\{\underbrace{q_1, \ldots, q_{k-2r}}_{\small{\hbox{part 2}}}
\}),
\end{equation}

\vspace{4mm}
\noindent
where 

\vspace{-2mm}
$$
\{g_1, g_2, \ldots, 
g_{2r-1}, g_{2r}, q_1, \ldots, q_{k-2r}\}=\{1, 2, \ldots, k\},
$$

\vspace{4mm}
\noindent
braces   
mean an unordered 
set, and pa\-ren\-the\-ses mean an ordered set.

We will say that (\ref{leto5007}) is a partition 
and consider the sum with respect to all possible
partitions

\vspace{1mm}
\begin{equation}
\label{leto5008}
\sum_{\stackrel{(\{\{g_1, g_2\}, \ldots, 
\{g_{2r-1}, g_{2r}\}\}, \{q_1, \ldots, q_{k-2r}\})}
{{}_{\{g_1, g_2, \ldots, 
g_{2r-1}, g_{2r}, q_1, \ldots, q_{k-2r}\}=\{1, 2, \ldots, k\}}}}
a_{g_1 g_2, \ldots, 
g_{2r-1} g_{2r}, q_1 \ldots q_{k-2r}},
\end{equation}

\vspace{5mm}
\noindent
where $a_{g_1 g_2, \ldots, 
g_{2r-1} g_{2r}, q_1 \ldots q_{k-2r}}\in \mathbb{R}.$

Below there are several examples of sums in the form (\ref{leto5008})

\vspace{2mm}
$$
\sum_{\stackrel{(\{g_1, g_2\})}{{}_{\{g_1, g_2\}=\{1, 2\}}}}
a_{g_1 g_2}=a_{12},
$$

\vspace{3mm}
$$
\sum_{\stackrel{(\{\{g_1, g_2\}, \{g_3, g_4\}\})}
{{}_{\{g_1, g_2, g_3, g_4\}=\{1, 2, 3, 4\}}}}
a_{g_1 g_2, g_3 g_4}=a_{12,34} + a_{13,24} + a_{23,14},
$$

\vspace{4mm}
$$
\sum_{\stackrel{(\{g_1, g_2\}, \{q_1, q_{2}\})}
{{}_{\{g_1, g_2, q_1, q_{2}\}=\{1, 2, 3, 4\}}}}
a_{g_1 g_2, q_1 q_{2}}=
$$

\vspace{1mm}
$$
=a_{12,34}+a_{13,24}+a_{14,23}
+a_{23,14}+a_{24,13}+a_{34,12},
$$

\vspace{4mm}
$$
\sum_{\stackrel{(\{g_1, g_2\}, \{q_1, q_{2}, q_3\})}
{{}_{\{g_1, g_2, q_1, q_{2}, q_3\}=\{1, 2, 3, 4, 5\}}}}
a_{g_1 g_2, q_1 q_{2}q_3}
=
$$

\vspace{1mm}
$$
=a_{12,345}+a_{13,245}+a_{14,235}
+a_{15,234}+a_{23,145}+a_{24,135}+
$$

\vspace{-2mm}
$$
+a_{25,134}+a_{34,125}+a_{35,124}+a_{45,123},
$$

\vspace{5mm}
$$
\sum_{\stackrel{(\{\{g_1, g_2\}, \{g_3, g_{4}\}\}, \{q_1\})}
{{}_{\{g_1, g_2, g_3, g_{4}, q_1\}=\{1, 2, 3, 4, 5\}}}}
a_{g_1 g_2, g_3 g_{4},q_1}
=
$$

\vspace{1mm}
$$
=
a_{12,34,5}+a_{13,24,5}+a_{14,23,5}+
a_{12,35,4}+a_{13,25,4}+a_{15,23,4}+
$$

\vspace{-2mm}
$$
+a_{12,54,3}+a_{15,24,3}+a_{14,25,3}+a_{15,34,2}+a_{13,54,2}+a_{14,53,2}+
$$

\vspace{-2mm}
$$
+
a_{52,34,1}+a_{53,24,1}+a_{54,23,1}.
$$

\vspace{6mm}

Now we can write (\ref{tyyy}) as

\vspace{1mm}

$$
J[\psi^{(k)}]_{T,t}=
\hbox{\vtop{\offinterlineskip\halign{
\hfil#\hfil\cr
{\rm l.i.m.}\cr
$\stackrel{}{{}_{p_1,\ldots,p_k\to \infty}}$\cr
}} }
\sum\limits_{j_1=0}^{p_1}\ldots
\sum\limits_{j_k=0}^{p_k}
C_{j_k\ldots j_1}\Biggl(
\prod_{l=1}^k\zeta_{j_l}^{(i_l)}+\sum\limits_{r=1}^{[k/2]}
(-1)^r \times
\Biggr.
$$

\vspace{3mm}
\begin{equation}
\label{leto6000hh}
\times
\sum_{\stackrel{(\{\{g_1, g_2\}, \ldots, 
\{g_{2r-1}, g_{2r}\}\}, \{q_1, \ldots, q_{k-2r}\})}
{{}_{\{g_1, g_2, \ldots, 
g_{2r-1}, g_{2r}, q_1, \ldots, q_{k-2r}\}=\{1, 2, \ldots, k\}}}}
\prod\limits_{s=1}^r
{\bf 1}_{\{i_{g_{{}_{2s-1}}}=~i_{g_{{}_{2s}}}\ne 0\}}
\Biggl.{\bf 1}_{\{j_{g_{{}_{2s-1}}}=~j_{g_{{}_{2s}}}\}}
\prod_{l=1}^{k-2r}\zeta_{j_{q_l}}^{(i_{q_l})}\Biggr),
\end{equation}

\vspace{5mm}
\noindent
where $[x]$ is an integer part of a real number $x,$
$\prod\limits_{\emptyset}
\stackrel{\sf def}{=}1,$ $\sum\limits_{\emptyset}
\stackrel{\sf def}{=}0;$
another notations are the same as in Theorem {\bf 1}.

\vspace{2mm}

In particular, from (\ref{leto6000hh}) for $k=5$ we obtain

\vspace{3mm}

$$
J[\psi^{(5)}]_{T,t}=
\hbox{\vtop{\offinterlineskip\halign{
\hfil#\hfil\cr
{\rm l.i.m.}\cr
$\stackrel{}{{}_{p_1,\ldots,p_5\to \infty}}$\cr
}} }\sum_{j_1=0}^{p_1}\ldots\sum_{j_5=0}^{p_5}
C_{j_5\ldots j_1}\Biggl(
\prod_{l=1}^5\zeta_{j_l}^{(i_l)}-\Biggr.
$$

\vspace{2mm}
$$
-
\sum\limits_{\stackrel{(\{g_1, g_2\}, \{q_1, q_{2}, q_3\})}
{{}_{\{g_1, g_2, q_{1}, q_{2}, q_3\}=\{1, 2, 3, 4, 5\}}}}
{\bf 1}_{\{i_{g_{{}_{1}}}=~i_{g_{{}_{2}}}\ne 0\}}
{\bf 1}_{\{j_{g_{{}_{1}}}=~j_{g_{{}_{2}}}\}}
\prod_{l=1}^{3}\zeta_{j_{q_l}}^{(i_{q_l})}+
$$

\vspace{2mm}
$$
+
\sum_{\stackrel{(\{\{g_1, g_2\}, 
\{g_{3}, g_{4}\}\}, \{q_1\})}
{{}_{\{g_1, g_2, g_{3}, g_{4}, q_1\}=\{1, 2, 3, 4, 5\}}}}
{\bf 1}_{\{i_{g_{{}_{1}}}=~i_{g_{{}_{2}}}\ne 0\}}
{\bf 1}_{\{j_{g_{{}_{1}}}=~j_{g_{{}_{2}}}\}}
\Biggl.{\bf 1}_{\{i_{g_{{}_{3}}}=~i_{g_{{}_{4}}}\ne 0\}}
{\bf 1}_{\{j_{g_{{}_{3}}}=~j_{g_{{}_{4}}}\}}
\zeta_{j_{q_1}}^{(i_{q_1})}\Biggr).
$$

\vspace{7mm}
\noindent
The last equality obviously agrees with
(\ref{a5}).

Let us consider the generalization of Theorem 1 for the case
of an arbitrary complete orthonormal systems  
of functions in the space $L_2([t,T])$ 
and $\psi_1(\tau),\ldots,\psi_k(\tau)\in L_2([t, T]).$

\vspace{2mm}

{\bf Theorem~2}\ \cite{20xxxxx} (Sect.~1.11), \cite{arxiv-1} (Sect.~15), \cite{new-2023a}.
{\it Suppose that
$\psi_1(\tau),\ldots,\psi_k(\tau)\in L_2([t, T])$ and
$\{\phi_j(x)\}_{j=0}^{\infty}$ is an arbitrary complete orthonormal system  
of functions in the space $L_2([t,T]).$
Then the following expansion

\vspace{2mm}
$$
J[\psi^{(k)}]_{T,t}=
\hbox{\vtop{\offinterlineskip\halign{
\hfil#\hfil\cr
{\rm l.i.m.}\cr
$\stackrel{}{{}_{p_1,\ldots,p_k\to \infty}}$\cr
}} }
\sum\limits_{j_1=0}^{p_1}\ldots
\sum\limits_{j_k=0}^{p_k}
C_{j_k\ldots j_1}\Biggl(
\prod_{l=1}^k\zeta_{j_l}^{(i_l)}+\sum\limits_{r=1}^{[k/2]}
(-1)^r \times
\Biggr.
$$

\vspace{4mm}
$$
\times
\sum_{\stackrel{(\{\{g_1, g_2\}, \ldots, 
\{g_{2r-1}, g_{2r}\}\}, \{q_1, \ldots, q_{k-2r}\})}
{{}_{\{g_1, g_2, \ldots, 
g_{2r-1}, g_{2r}, q_1, \ldots, q_{k-2r}\}=\{1, 2, \ldots, k\}}}}
\prod\limits_{s=1}^r
{\bf 1}_{\{i_{g_{{}_{2s-1}}}=~i_{g_{{}_{2s}}}\ne 0\}}
\Biggl.{\bf 1}_{\{j_{g_{{}_{2s-1}}}=~j_{g_{{}_{2s}}}\}}
\prod_{l=1}^{k-2r}\zeta_{j_{q_l}}^{(i_{q_l})}\Biggr)
$$

\vspace{5mm}
\noindent
con\-verg\-ing in the mean-square sense is valid,
where $[x]$ is an integer part of a real number $x,$
$\prod\limits_{\emptyset}
\stackrel{\sf def}{=}1,$ $\sum\limits_{\emptyset}
\stackrel{\sf def}{=}0;$
another notations are the same as in Theorem~{\rm 1}.}

\vspace{2mm}

It should be noted that an analogue of Theorem 2 was considered 
in \cite{Rybakov1000} using Hermite polynomials. 
Note that we use another notations 
\cite{20xxxxx} (Sect.~1.11), \cite{arxiv-1} (Sect.~15), \cite{new-2023a}
in comparison with \cite{Rybakov1000}.
Moreover, the proof
from \cite{Rybakov1000} is different from the proof given in 
\cite{20xxxxx} (Sect.~1.11), \cite{arxiv-1} (Sect.~15), \cite{new-2023a}.
The results of \cite{Rybakov1000}, as well as 
the results of \cite{20xxxxx} (Sect.~1.11), \cite{arxiv-1} (Sect.~15), \cite{new-2023a}
are based on our idea 
\cite{7} (2006) on the expansion of the kernel (\ref{ppp}) (or $\Phi(t_1,\ldots,t_k)\in L_2([t,T]^k)$)
into a generalized multiple Fourier series 
(see \cite{7}, Chapter~5, Theorem~5.1, pp.~235-245 
or \cite{20xxxxx}, Chapter~1 for details).

\vspace{5mm}

\section{Theorem on Expansion of Iterated 
Stratonovich Stochastic Integrals of Second Multiplicity. Some Old Results}

\vspace{5mm}

In a number of works of the author 
\cite{12}-\cite{16}, \cite{19}-\cite{21},
\cite{23}, \cite{25}, \cite{27}, 
\cite{arxiv-7}, \cite{arxiv-8ee}, \cite{100-000-1},
\cite{100-000-5}
Theorems 1, 2 have been adapted for iterated Stratonovich stochastic integrals
(\ref{605}) of multiplicities 2 to 6 (also see the case of multiplicity
$k$ $(k\in\mathbb{N})$ in \cite{20xxxxx} (Sect.~2.10), \cite{21}, \cite{25}, \cite{arxiv-4}).
For example, we can formulate the following theorem for iterated 
Stratonovich stochastic integrals of second multiplicity.

\vspace{2mm}

{\bf Theorem 3} \cite{12}-\cite{16}, \cite{19}-\cite{21},
\cite{arxiv-7}, \cite{arxiv-8ee}, \cite{100-000-1},
\cite{100-000-5}.
{\it Suppose that 
$\{\phi_j(x)\}_{j=0}^{\infty}$ is a complete or\-tho\-nor\-mal system of 
Legendre polynomials or trigonometric functions in the space $L_2([t, T]).$
At the same time $\psi_2(\tau)$ is a continuously differentiable 
function on $[t, T]$ and $\psi_1(\tau)$ is twice
continuously differentiable functions on $[t, T]$. Then

\vspace{1mm}
\begin{equation}
\label{a}
J^{*}[\psi^{(2)}]_{T,t}=
\hbox{\vtop{\offinterlineskip\halign{
\hfil#\hfil\cr
{\rm l.i.m.}\cr
$\stackrel{}{{}_{p_1,p_2\to \infty}}$\cr
}} }\sum_{j_1=0}^{p_1}\sum_{j_2=0}^{p_2}
C_{j_2j_1}\zeta_{j_1}^{(i_1)}\zeta_{j_2}^{(i_2)}\ \ \ (i_1,i_2=1,\ldots,m),
\end{equation}

\vspace{4mm}
\noindent
where $J^{*}[\psi^{(2)}]_{T,t}$ is defined by {\rm (\ref{605});} 
another notations are the
same as in Theorem {\rm 1.}
}

\vspace{4mm}

Note that the proof of Theorem 3 is based
on the proof of the following equality

\vspace{1mm}
\begin{equation}
\label{5t11}
\frac{1}{2}
\int\limits_t^T\psi_1(t_1)\psi_2(t_1)dt_1
=\sum_{j_1=0}^{\infty}
C_{j_1j_1},
\end{equation}

\vspace{4mm}
\noindent
where $C_{j_1j_1}$ is defined by (\ref{ppppa})
for $k=2$ and $j_1=j_2;$
$\{\phi_j(x)\}_{j=0}^{\infty}$ is a complete orthonormal system of 
Legendre polynomials or trigonometric functions in the space $L_2([t, T]).$

According to the standard relation between 
Ito and Stratonovich stochastic integrals, we can write w.~p.~1
(with probability 1)

\vspace{-1mm}
\begin{equation}
\label{oop51}
J^{*}[\psi^{(2)}]_{T,t}=
J[\psi^{(2)}]_{T,t}+
\frac{1}{2}{\bf 1}_{\{i_1=i_2\ne 0\}}
\int\limits_t^T\psi_1(t_1)\psi_2(t_1)dt_1,
\end{equation}

\vspace{3mm}
\noindent
where we assume that the functions $\psi_1(\tau), \psi_2(\tau)$ are continuous
at the interval $[t, T].$
This condition is related to the definition of the Stratonovich stochastic integral 
that we use \cite{KlPl2} (also see Sect.~2.1.1 \cite{20xxxxx}).

From the other hand according to (\ref{leto5001}), we obtain

\vspace{2mm}
$$
J[\psi^{(2)}]_{T,t}=
\hbox{\vtop{\offinterlineskip\halign{
\hfil#\hfil\cr
{\rm l.i.m.}\cr
$\stackrel{}{{}_{p_1,p_2\to \infty}}$\cr
}} }\sum_{j_1=0}^{p_1}\sum_{j_2=0}^{p_2}
C_{j_2j_1}\Biggl(\zeta_{j_1}^{(i_1)}\zeta_{j_2}^{(i_2)}
-{\bf 1}_{\{i_1=i_2\ne 0\}}
{\bf 1}_{\{j_1=j_2\}}\Biggr)=
$$

\vspace{4mm}
\begin{equation}
\label{yes2001}
=\hbox{\vtop{\offinterlineskip\halign{
\hfil#\hfil\cr
{\rm l.i.m.}\cr
$\stackrel{}{{}_{p_1,p_2\to \infty}}$\cr
}} }\sum_{j_1=0}^{p_1}\sum_{j_2=0}^{p_2}
C_{j_2j_1}\zeta_{j_1}^{(i_1)}\zeta_{j_2}^{(i_2)}
-{\bf 1}_{\{i_1=i_2\ne 0\}}\sum_{j_1=0}^{\infty}
C_{j_1j_1}.
\end{equation}

\vspace{5mm}

From (\ref{5t11})--(\ref{yes2001}) we get (\ref{a}). Note that
the existence of the limit on the right-hand side of (\ref{5t11})
will be proved below (see Lemma 2 and Theorem 7).

The proof of Theorem 3
\cite{12}-\cite{16}, \cite{19}-\cite{21},
\cite{23}, \cite{25}, \cite{27}, 
\cite{arxiv-7}, \cite{arxiv-8ee}, \cite{100-000-1},
\cite{100-000-5}
is based on double (iterated ) Fourier--Legendre
series and analogous trigonometric Fourier series.
This proof uses double integration by parts, which leads to 
the requirement of double continuous differentiability of the function
$\psi_1(\tau)$ 
at the interval $[t, T]$.

In this article, we formulate and prove the analogue
of Theorem 3 (Theorem 6, see below) 
but under the weakened conditions:
the functions
$\psi_1(\tau),$  $\psi_2(\tau)$
are assumed to be continuously differentiable only one time
at the interval $[t, T]$.
At that we will use 
double Fourier--Legendre series
and double trigonometric Fourier series summarized by Pringsheim method
for the proof of Theorem 6 (see below).

In Sect.~5 (see Theorem~7), we generalize the equality (\ref{5t11}) 
to the case of an arbitrary 
complete orthonormal system of functions 
in the space $L_2([t, T])$ and $\psi_1(\tau),\psi_2(\tau)\in L_2([t,T]).$

\vspace{5mm}

\section{Proof of the Equality (\ref{5t11}). 
The Case of Legendre Polynomials and Trigonometric Functions
as well as Continuously Differentiable Functions $\psi_1(\tau),$ $\psi_2(\tau)$}

\vspace{5mm}

Let $P_j(x)$ $(j=0, 1, 2,\ldots )$ be the Legendre polynomial.
Consider the function $f(x,y)$ defined for $(x,y)\in [-1,1]^2.$

Consider the double Fourier--Legendre series summarized by
Pringsheim method and corresponding 
to the function 
$f(x,y)$

\vspace{1mm}
$$
\lim_{n,m\to\infty}
\sum_{j=0}^n\sum_{i=0}^m \frac{\sqrt{(2j+1)(2i+1)}}{2}{C}_{ij}^{*}
P_i(x)P_j(y)
\stackrel{\sf{def}}{=}
$$

\vspace{3mm}
\begin{equation}
\label{555}
\stackrel{\sf{def}}{=}
\sum_{i,j=0}^{\infty}\frac{\sqrt{(2j+1)(2i+1)}}{2}{C}_{ij}^{*}
P_i(x)P_j(y),
\end{equation}

\vspace{3mm}
\noindent
where

\vspace{-1mm}
\begin{equation}
\label{777}
{C}_{ij}^{*}=\frac{\sqrt{(2j+1)(2i+1)}}{2}\int\limits_{[-1,1]^2}
f(x,y)P_i(x)P_j(y)dxdy.
\end{equation}

\vspace{4mm}

Let us consider the generalization for the case of two variables
\cite{star}
of the theorem on equiconvergence for the Fourier--Legendre
series
\cite{suet}.

\vspace{2mm}

{\bf Theorem 4}\ \cite{star}.\ {\it Let $f(x,y)\in L_2([-1,1]^2)$ and 
the function

\vspace{1mm}
$$
f(x,y)(1-x^2)^{-1/4}(1-y^2)^{-1/4}
$$

\vspace{4mm}
\noindent
is integrable on the square
$[-1,1]^2.$  Moreover, let

\vspace{1mm}
$$
|f(x,y)-f(u,v)|\le G(y)|x-u|+H(x)|y-v|,
$$

\vspace{4mm}
\noindent
where $G(y), H(x)$ are bounded functions on the square $[-1,1]^2$.
Then for all $(x,y)\in(-1,1)^2$ the following equality
is satisfied 

\vspace{4mm}
$$
\lim\limits_{n,m\to\infty}
\Biggl(\sum_{j=0}^n\sum_{i=0}^m \frac{\sqrt{(2j+1)(2i+1)}}{2}{C}_{ij}^{*}
P_i(x)P_j(y)-\Biggr.
$$

\begin{equation}
\label{888}
-\Biggl.
(1-x^2)^{-1/4}(1-y^2)^{-1/4}S_{nm}({\rm arccos}x,{\rm arccos}y,F)
\Biggr)=0,
\end{equation}

\vspace{4mm}
\noindent
where $S_{nm}(\theta,\varphi,F)$ is a partial sum of
the double trigonometric Fourier series of the
auxiliary function

\vspace{2mm}
$$
F(\theta,\varphi)=\sqrt{|{\rm sin}\theta|}\sqrt{|{\rm sin}\varphi|}
f({\rm cos}\theta,{\rm cos}\varphi),\ \ \
\theta,\varphi\in[0, \pi],
$$

\vspace{5mm}
\noindent
the Fourier coefficient ${C}_{ij}^{*}$ has the form {\rm (\ref{777}).}
At that, the convergence in {\rm (\ref{888})} is uniform on the rectangle

\vspace{-1mm}
$$
[-1+\varepsilon, 1-\varepsilon]\times[-1+\delta,1-\delta]\ \
\ \hbox{for any}\ \ \ \varepsilon, \delta>0.
$$ 
}

\vspace{2mm}

From Theorem 4, for example, follows that  
for all $(x,y)\in(-1,1)^2$ the following equality is
fulfilled 

\vspace{1mm}
\begin{equation}
\label{8888}
\lim\limits_{n,m\to\infty}
\Biggl(\sum_{j=0}^n\sum_{i=0}^m \frac{\sqrt{(2j+1)(2i+1)}}{2}{C}_{ij}^{*}
P_i(x)P_j(y)-f(x,y)\Biggr)=0
\end{equation}

\vspace{5mm}
\noindent
and convergence in {\rm (\ref{8888})} is uniform on the rectangle

$$
[-1+\varepsilon, 1-\varepsilon]\times[-1+\delta,1-\delta]\ \ \
\hbox{for any}\ \ \ \varepsilon, \delta>0
$$ 

\vspace{4mm}
\noindent
if the corresponding conditions of convergence 
of the double trigonometric Fourier series of the auxiliary 
function 

\vspace{-1mm}
\begin{equation}
\label{999}
g(x,y)=f(x,y)(1-x^2)^{1/4}(1-y^2)^{1/4}
\end{equation}

\vspace{3mm}
\noindent
are satisfied.

Note that Theorem 4 does not imply any
conclusions on the behavior of the double Fourier--Legendre series 
on the boundary of
the square 
$[-1,1]^2.$

\vspace{2mm}

{\it For each $\delta>0$ let us call the exact upper edge 
of the difference $|f({\bf t}')-f({\bf t}'')|$ 
in the set
of all points ${\bf t}'$, ${\bf t}''$ {\rm (}which 
belong 
to the domain $D${\rm )}
as the module of 
continuity of the function
$f({\bf t})$ $({\bf t}=(t_1,\ldots,t_k))$ in the 
$k$-dimentional domain
$D$ $(k\ge 1)$ if
the distance $\rho({\bf t}',{\bf t}'')$ between 
${\bf t}$ and ${\bf t}''$
satisfies the condition
$\rho({\bf t}',{\bf t}'')<\delta.$}

\vspace{2mm}

{\it We will say that the function 
of $k$ $(k\ge 1)$ variables  
$f({\bf t})$ $({\bf t}=(t_1,\ldots,t_k))$
belongs 
to the H\"{o}lder class with 
the parameter $\alpha\in (0, 1]$ $(f({\bf t})\in C^{\alpha}(D))$ 
in the domain $D$ 
if the module of 
continuity of the function
$f({\bf t})$ $({\bf t}=(t_1,\ldots,t_k))$ 
in the domain $D$ has orders $o(\delta^{\alpha})$ $(\alpha \in (0, 1))$
and $O(\delta)$ $(\alpha=1)$.
}

\vspace{2mm}

In 1967, Zhizhiashvili L.V.
proved
that the rectangular sums of multiple trigonometric Fourier series 
of the function of $k$ variables 
in the hypercube $[t,T]^k$ 
converge uniformly in the hypercube $[t,T]^k$
to this function if it 
belongs
to the class $C^{\alpha}([t,T]^k),$ $\alpha>0$ (definition
of the Holder class with the parameter $\alpha>0$ can be found in 
the well-known mathematical analysis tutorials; see, for
example, \cite{IP}).
More precisely, the following theorem is correct.

\vspace{2mm}

{\bf Theorem 5}\ \cite{IP}.\ {\it If the function
$f(x_1,\ldots,x_n)$ is periodic with period $2\pi$ with respect to each
variable and belongs in $\mathbb{R}^n$ to the Holder class $C^{\alpha}$
for any $\alpha>0,$ then the rectangular partial
sums of the multiple trigonometric Fourier series of the function
$f(x_1,\ldots,x_n)$ converge to this function uniformly in 
$\mathbb{R}^n$.}

\vspace{2mm}

{\bf Lemma 1}.\ {\it 
Let the function $f(x,y)$ satisfies
to the following condition 

\vspace{1mm}
$$
|f(x,y)-f(x_1,y_1)|\le C_1|x-x_1|+C_2|y-y_1|,
$$

\vspace{5mm}
\noindent
where $C_1, C_2<\infty$ and $(x,y),$ $(x_1,y_1)\in [-1,1]^2.$
Then the following inequality is fulfilled

\begin{equation}
\label{2222}
|g(x,y)-g(x_1,y_1)|\le K\rho^{1/4},
\end{equation}

\vspace{4mm}
\noindent
where $g(x,y)$ 
has the form
{\rm(\ref{999})}{\rm,}\ 

$$
\rho=\sqrt{(x-x_1)^2+(y-y_1)^2},
$$ 

\vspace{5mm}
\noindent
$(x,y)$ and $(x_1,y_1)\in [-1,1]^2,$ $K<\infty.$}
                                  
\vspace{2mm}

{\bf Proof.}\ First, we assume that $x\ne x_1,$ $y\ne y_1.$ In this case
we have

\vspace{2mm}
$$
|g(x,y)-g(x_1,y_1)|=
|(1-x^2)^{1/4}(1-y^2)^{1/4}(f(x,y)-f(x_1,y_1))+
$$

\vspace{1mm}
$$
+
f(x_1,y_1)((1-x^2)^{1/4}(1-y^2)^{1/4}
-(1-x_1^2)^{1/4}(1-y_1^2)^{1/4})|\le C_1|x-x_1|+C_2|y-y_1|+
$$

\vspace{1mm}
\begin{equation}
\label{6666}
+C_3|(1-x^2)^{1/4}(1-y^2)^{1/4}
-(1-x_1^2)^{1/4}(1-y_1^2)^{1/4}|,\ \ \ C_3<\infty.
\end{equation}

\vspace{7mm}

Moreover,

\vspace{-1mm}
$$
|(1-x^2)^{1/4}(1-y^2)^{1/4}
-(1-x_1^2)^{1/4}(1-y_1^2)^{1/4}|=
$$

\vspace{1mm}
$$
=|(1-x^2)^{1/4}((1-y^2)^{1/4}-(1-y_1^2)^{1/4})
+(1-y_1^2)^{1/4}((1-x^2)^{1/4}-(1-x_1^2)^{1/4})|\le
$$

\vspace{1mm}
\begin{equation}
\label{5555}
\le |(1-y^2)^{1/4}-(1-y_1^2)^{1/4}|+
|(1-x^2)^{1/4}-(1-x_1^2)^{1/4}|,
\end{equation}

\vspace{5mm}

$$
|(1-x^2)^{1/4}-(1-x_1^2)^{1/4}|=
$$

$$
=|((1-x)^{1/4}-(1-x_1)^{1/4})(1+x)^{1/4}+
(1-x_1)^{1/4}((1+x)^{1/4}-(1+x_1)^{1/4})|\le
$$

\begin{equation}
\label{4444}
\le K_1(|(1-x)^{1/4}-(1-x_1)^{1/4}|+
|(1+x)^{1/4}-(1+x_1)^{1/4}|),\ \ \ K_1<\infty.
\end{equation}

\vspace{7mm}

It is not difficult to see that

$$
|(1 \pm x)^{1/4}-(1 \pm x_1)^{1/4}|=
$$

\vspace{3mm}
$$
=
\frac{|(1\pm x) - (1\pm x_1)|}
{((1\pm x)^{1/2}+(1 \pm x_1)^{1/2})((1\pm x)^{1/4}+(1 \pm x_1)^{1/4})}=
$$

\vspace{3mm}
\begin{equation}
\label{3333}
=|x_1-x|^{1/4} \frac{|x_1-x|^{1/2}}{(1\pm x)^{1/2}+(1 \pm x_1)^{1/2}}
\cdot \frac{|x_1-x|^{1/4}}{(1\pm x)^{1/4}+(1 \pm x_1)^{1/4}}
\le |x_1-x|^{1/4}.
\end{equation}

\vspace{9mm}

The last inequality follows from the obvious inequalities

\vspace{2mm}
$$
\frac{|x_1-x|^{1/2}}{(1\pm x)^{1/2}+(1 \pm x_1)^{1/2}}\le 1,
$$

\vspace{3mm}
$$
\frac{|x_1-x|^{1/4}}{(1\pm x)^{1/4}+(1 \pm x_1)^{1/4}}\le 1.
$$

\vspace{6mm}

From (\ref{6666})--(\ref{3333}) we obtain

$$
|g(x,y)-g(x_1,y_1)|\le
C_1|x-x_1|+C_2|y-y_1|+C_4(|x_1-x|^{1/4}+|y_1-y|^{1/4})\le
$$

\vspace{1mm}
$$
\le 
C_5 \rho + C_6 \rho^{1/4} \le K \rho^{1/4},
$$

\vspace{4mm}
\noindent
where $C_5, C_6, K <\infty.$

The cases $x=x_1, y\ne y_1$ and $x\ne x_1, y=y_1$ 
can be considered analogously to the case
$x\ne x_1, y\ne y_1$. At that, the consideration 
begins from the inequalities 

\vspace{1mm}
$$
|g(x,y)-g(x_1,y_1)|\le
K_2 |(1-y^2)^{1/4}f(x,y)-
(1-y_1^2)^{1/4}f(x_1,y_1)|
$$

\vspace{5mm}
\noindent
($x=x_1,\ y\ne y_1$) and

\vspace{1mm}
$$
|g(x,y)-g(x_1,y_1)|\le  
K_2 |(1-x^2)^{1/4}f(x,y)-
(1-x_1^2)^{1/4}f(x_1,y_1)|
$$

\vspace{5mm}
\noindent
($x\ne x_1,\ y=y_1$),
where $K_2<\infty$. Lemma 1 is proved. 

\vspace{2mm}

Lemma 1 and Theorem 5 imply
that rectangular partial sums of the double trigonometric
Fourier series of the function $g (x,y)$ (in the case 
of periodic continuation of the function 
$g (x,y)$)
converge uniformly in 
the square $[-1,1]^2$
to the function $g (x,y)$. This means that
the equality (\ref{8888}) holds.

\vspace{2mm}

{\bf Theorem 6}\ \cite{20xxxxx}-\cite{12aa-afterxxx}, \cite{30a}, \cite{arxiv-8}. {\it 
Suppose that 
$\{\phi_j(x)\}_{j=0}^{\infty}$ is a complete orthonormal system of 
Legendre polynomials or trigonometric functions in the space $L_2([t, T]).$
Moreover, $\psi_1(\tau),$ $\psi_2(\tau)$ are 
continuously differentiable functions on $[t, T]$. Then 
for the iterated Stratonovich stochastic integral

$$
J^{*}[\psi^{(2)}]_{T,t}={\int\limits_t^{*}}^T\psi_2(t_2)
{\int\limits_t^{*}}^{t_2}\psi_1(t_1)d{\bf f}_{t_1}^{(i_1)}
d{\bf f}_{t_2}^{(i_2)}\ \ \ (i_1, i_2=1,\ldots,m)
$$

\vspace{3mm}
\noindent
the following expansion

\begin{equation}
\label{jes}
J^{*}[\psi^{(2)}]_{T,t}=\hbox{\vtop{\offinterlineskip\halign{
\hfil#\hfil\cr
{\rm l.i.m.}\cr
$\stackrel{}{{}_{p_1,p_2\to \infty}}$\cr
}} }\sum_{j_1=0}^{p_1}\sum_{j_2=0}^{p_2}
C_{j_2j_1}\zeta_{j_1}^{(i_1)}\zeta_{j_2}^{(i_2)}
\end{equation}

\vspace{4mm}
\noindent
that converges in the mean-square
sence is valid, where the notations are the same as in Theorem {\rm 3.}
}

\vspace{2mm}

{\bf Proof.}\ Let us prove the equality

$$
\frac{1}{2}
\int\limits_t^T\psi_1(t_1)\psi_2(t_1)dt_1
=\sum_{j_1=0}^{\infty}
C_{j_1j_1},
$$

\vspace{3mm}
\noindent
where $C_{j_1j_1}$ is defined by the formula (\ref{ppppa})
for $k=2$ and $j_1=j_2.$ At that
$\{\phi_j(x)\}_{j=0}^{\infty}$ is a complete orthonormal system of 
Legendre polynomials or trigonometric functions in the space $L_2([t, T]).$

Consider the auxiliary function

$$
K'(t_1,t_2)=\left\{
\begin{matrix}
\psi_2(t_1)\psi_1(t_2),\ \ t_1\ge t_2\cr\cr\cr
\psi_1(t_1)\psi_2(t_2),\ \ t_1\le t_2
\end{matrix}
\right.,\ \ \ t_1,t_2\in[t,T]
$$

\vspace{3mm}
\noindent
and prove that 

\vspace{-2mm}
\begin{equation}
\label{9090}
|K'(t_1,t_2)-K'(t_1^{*},t_2^{*})|\le L(|t_1-t_1^{*}|+|t_2-t_2^{*}|),
\end{equation}

\vspace{4mm}
\noindent
where $L<\infty,$ and $(t_1,t_2)$, $(t_1,t_2)\in [t, T]^2.$

By the Lagrange formula for
the functions $\psi_1(t_1^{*}),$ $\psi_2(t_1^{*})$ at the interval
$[{\rm min}\{t_1, t_1^{*}\}, {\rm max}\{t_1, t_1^{*}\}]$ 
and for the functions 
$\psi_1(t_2^{*}),$ $\psi_2(t_2^{*})$ at the interval
$[{\rm min}\{t_2, t_2^{*}\}, {\rm max}\{t_2, t_2^{*}\}]$ 
we obtain

\vspace{2mm}
$$
|K'(t_1,t_2)-K'(t_1^{*},t_2^{*})|\le 
\left|\ \left\{
\begin{matrix}
\psi_2(t_1)\psi_1(t_2),\ \ t_1\ge t_2\cr\cr\cr
\psi_1(t_1)\psi_2(t_2),\ \ t_1\le t_2
\end{matrix}\right.
-\left\{
\begin{matrix}
\psi_2(t_1)\psi_1(t_2),\ t_1^{*}\ge t_2^{*}\cr\cr\cr
\psi_1(t_1)\psi_2(t_2),\ \ t_1^{*}\le t_2^{*}
\end{matrix}
\right.\ \right|+
$$

\vspace{3mm}
\begin{equation}
\label{uuuu1}
+
L_1|t_1-t_1^{*}|+L_2|t_2-t_2^{*}|,\ \ \ L_1, L_2<\infty.
\end{equation}

\vspace{7mm}

We have

$$
\left\{
\begin{matrix}
\psi_2(t_1)\psi_1(t_2),\ \ t_1\ge t_2\cr\cr\cr
\psi_1(t_1)\psi_2(t_2),\ \ t_1\le t_2
\end{matrix}
\right.
-\left\{
\begin{matrix}
\psi_2(t_1)\psi_1(t_2),\ \ t_1^{*}\ge t_2^{*}\cr\cr\cr
\psi_1(t_1)\psi_2(t_2),\ \ t_1^{*}\le t_2^{*}
\end{matrix}\right.
=
$$

\vspace{6mm}
\begin{equation}
\label{uuuu2}
=\left\{
\begin{matrix}
0,\ \ \ t_1\ge t_2,\  t_1^{*}\ge t_2^{*}\ \ \ \hbox{or}\ \ \ 
t_1\le t_2,\ t_1^{*}\le t_2^{*}\cr\cr\cr
\psi_2(t_1)\psi_1(t_2)-\psi_1(t_1)\psi_2(t_2),\ \ \
t_1\ge t_2,\ t_1^{*}\le t_2^{*}\cr\cr\cr
\psi_1(t_1)\psi_2(t_2)-\psi_2(t_1)\psi_1(t_2),\ \ \
t_1\le t_2,\ t_1^{*}\ge t_2^{*}
\end{matrix}.\right.
\end{equation}

\vspace{7mm}

By the Lagrange formula for the functions
$\psi_1(t_2),$ $\psi_2(t_2)$ at the interval
$[{\rm min}\{t_1, t_2\}, {\rm max}\{t_1, t_2\}]$ 
we get the estimate

\vspace{1mm}
$$
\left|\ \left\{
\begin{matrix}
\psi_2(t_1)\psi_1(t_2),\ \ t_1\ge t_2\cr\cr\cr
\psi_1(t_1)\psi_2(t_2),\ \ t_1\le t_2
\end{matrix}\right.
-\left\{
\begin{matrix}\psi_2(t_1)\psi_1(t_2),\ \ t_1^{*}\ge t_2^{*}\cr\cr\cr
\psi_1(t_1)\psi_2(t_2),\ \ t_1^{*}\le t_2^{*}
\end{matrix}
\right.\ \right|\le
$$

\vspace{5mm}
\begin{equation}
\label{ppp1}
\le L_3 |t_2-t_1|
\left\{
\begin{matrix}
0,\ \ \ t_1\ge t_2,\ t_1^{*}\ge t_2^{*}\ \ \ \hbox{or}\ \ \
t_1\le t_2,\ t_1^{*}\le t_2^{*}\cr\cr\cr
1,\ \ \ t_1\le t_2,\ t_1^{*}\ge t_2^{*}\ \ \ \hbox{or}\ \ \
t_1\ge t_2,\ t_1^{*}\le t_2^{*}
\end{matrix},\right.\ \ \ L_3<\infty. 
\end{equation}

\vspace{7mm}

Let us show that if $t_1\le t_2,\ t_1^{*}\ge t_2^{*}$ or
$t_1\ge t_2,\ t_1^{*}\le t_2^{*}$, then the following 
inequality is satisfied  

\vspace{1mm}
\begin{equation}
\label{r5r5}
|t_2-t_1|\le |t_1^{*}-t_1|+|t_2^{*}-t_2|.
\end{equation}

\vspace{4mm}

First, consider the case $t_1\ge t_2,\ t_1^{*}\le t_2^{*}$.
For this case

\vspace{1mm}
$$
t_2+(t_1^{*}-t_2^{*})\le t_2\le t_1.
$$

\vspace{4mm}

Then 
$$
(t_1^{*}-t_1)-(t_2^{*}-t_2)\le t_2-t_1\le 0
$$ 

\vspace{3mm}
\noindent
and (\ref{r5r5}) is satisfied.

For the case $t_1\le t_2,\ t_1^{*}\ge t_2^{*}$ we have

\vspace{1mm}
$$
t_1+(t_2^{*}-t_1^{*})\le t_1\le t_2.
$$

\vspace{4mm}

Then 
$$
(t_1-t_1^{*})-(t_2-t_2^{*})\le t_1-t_2\le 0
$$ 

\vspace{3mm}
\noindent
and (\ref{r5r5}) is also satisfied.

From (\ref{ppp1}) and (\ref{r5r5}) we obtain

\vspace{3mm}
$$
\left|\ \left\{
\begin{matrix}
\psi_2(t_1)\psi_1(t_2),\ \ t_1\ge t_2\cr\cr\cr
\psi_1(t_1)\psi_2(t_2),\ \ t_1\le t_2
\end{matrix}\right.
-\left\{
\begin{matrix}
\psi_2(t_1)\psi_1(t_2),\ \ t_1^{*}\ge t_2^{*}\cr\cr\cr
\psi_1(t_1)\psi_2(t_2),\ \ t_1^{*}\le t_2^{*}
\end{matrix}
\right.\ \right|\le
$$

\vspace{5mm}
$$
\le L_3 (|t_1^{*}-t_1|+|t_2^{*}-t_2|)
\left\{
\begin{matrix}
0,\ \ \ t_1\ge t_2,\ t_1^{*}\ge t_2^{*}\ \ \ \hbox{or}\ \ \ 
t_1\le t_2,\ t_1^{*}\le t_2^{*}\cr\cr\cr
1,\ \ \ t_1\le t_2,\ t_1^{*}\ge t_2^{*}\ \ \ \hbox{or}\ \ \ 
t_1\ge t_2,\ t_1^{*}\le t_2^{*}
\end{matrix}\ \right. \le
$$

\vspace{5mm}
$$
\le
L_3 (|t_1^{*}-t_1|+|t_2^{*}-t_2|)
\left\{
\begin{matrix}
1,\ \ \ t_1\ge t_2,\ t_1^{*}\ge t_2^{*}\ \ \ \hbox{or}\ \ \ 
t_1\le t_2,\ t_1^{*}\le t_2^{*}\cr\cr\cr
1,\ \ \ t_1\le t_2,\ t_1^{*}\ge t_2^{*}\ \ \ \hbox{or}\ \ \ 
t_1\ge t_2,\ t_1^{*}\le t_2^{*}
\end{matrix}=\right. 
$$

\vspace{3mm}
\begin{equation}
\label{t2t2}
=L_3 (|t_1^{*}-t_1|+|t_2^{*}-t_2|).
\end{equation}

\vspace{7mm}

From (\ref{uuuu1}) and (\ref{t2t2}) we get (\ref{9090}).
Let

\vspace{1mm}
$$
t_1=\frac{T-t}{2}x+\frac{T+t}{2},\ \ \
t_2=\frac{T-t}{2}y+\frac{T+t}{2},
$$

\vspace{4mm}
\noindent
where $x, y \in [-1,1].$ 
Then

$$
K'(t_1,t_2)\equiv K^{*}(x,y)
=\left\{
\begin{matrix}
\psi_2\left(h(x)\right)
\psi_1\left(h(y)\right),\ \ x\ge y\cr\cr\cr
\psi_1\left(h(x)\right)
\psi_2\left(h(y)\right),\ \ x\le y
\end{matrix}\right.,
$$

\vspace{5mm}
\noindent
where $x, y\in[-1,1]$ and

\vspace{-3mm}
\begin{equation}
\label{ter1}
h(x)=\frac{T-t}{2}x+\frac{T+t}{2}.
\end{equation}

\vspace{5mm}
                               
Inequality (\ref{9090}) can be written in the form 

\begin{equation}
\label{4343}
|K^{*}(x,y)- K^{*}(x^{*},y^{*})|\le L^{*}(|x-x^{*}|+|y-y^{*}|),
\end{equation}

\vspace{4mm}
\noindent
where $L^{*}<\infty$ and $(x,y)$, $(x^{*},y^{*})\in [-1, 1]^2.$

Thus, the function $K^{*} (x,y)$ satisfies the conditions of 
Lemma 1 and hence for the function 

\vspace{1mm}
$$
K^{*}(x,y)(1-x^2)^{1/4}(1-y^2)^{1/4}
$$

\vspace{4mm}
\noindent
the inequality (\ref{2222}) is fulfilled.

Due to the continuous differentiability of the functions 
$\psi_1\left(h(x)\right)$ and
$\psi_2\left(h(x)\right)$ at the interval
$[-1,1]$ we have
$K^{*}(x,y)\in L_2([-1,1]^2)$. In addition

\vspace{2mm}
$$
\int\limits_{[-1,1]^2}
\frac{K^{*}(x,y)dxdy}{(1-x^2)^{1/4}(1-y^2)^{1/4}}\le
C \left(~\int\limits_{-1}^{1}\frac{1}{(1-x^2)^{1/4}}
\int\limits_{-1}^{x}\frac{1}{(1-y^2)^{1/4}}dydx+\right.
$$

\vspace{4mm}
$$
\left.+
\int\limits_{-1}^{1}\frac{1}{(1-x^2)^{1/4}}
\int\limits_{x}^{1}\frac{1}{(1-y^2)^{1/4}}dydx\right)<\infty,\ \ \ C<\infty.
$$

\vspace{6mm}

Thus, the conditions of Theorem 4 are fulfilled
for the function $K^{*} (x,y).$ 
Note that the mentioned properties of the function
$K^{*}(x,y),$ $x,y \in [-1,1]$ also correct
for the function $K'(t_1,t_2),$ $t_1, t_2\in [t, T].$

\vspace{2mm}

{\bf Remark 1.} {\it On the basis of {\rm(\ref{9090})} it can be 
argued that the 
function $K'(t_1, t_2)$ belongs to the Holder class 
with parameter $1$ in $[t, T]^2.$ Hence by Theorem $5$ this 
function can be expanded into the uniformly convergent double 
trigonometric Fourier series in the square $[t, T]^2,$ 
which summarized by Pringsheim method.
However, the expansions of iterated stochastic integrals 
obtained by using the system of
Legendre polynomials are essentially simpler than their 
analogues obtained by using the trigonometric system of functions
{\rm (}see Sect.~{\rm 7)}.}

\vspace{2mm}

Let us expand the function $K'(t_1,t_2)$ into a
double Fourier--Legendre series or double trigonometric Fourier series  
in the square $[t, T]^2$. This 
series is summable 
by the method of rectangular sums (Pringsheim method), i.e.

\vspace{1mm}
$$
K'(t_1,t_2)=
\lim_{n_1,n_2\to\infty}
\sum_{j_1=0}^{n_1}\sum_{j_2=0}^{n_2}
\int\limits_t^T\int\limits_t^TK'(t_1,t_2)
\phi_{j_1}(t_1)\phi_{j_2}(t_2)dt_1 dt_2\cdot
\phi_{j_1}(t_1)\phi_{j_2}(t_2)=
$$

\vspace{2mm}
$$
=\lim_{n_1,n_2\to\infty}
\sum_{j_1=0}^{n_1}\sum_{j_2=0}^{n_2}\left(
\int\limits_t^T\psi_2(t_2)\phi_{j_2}(t_2)
\int\limits_t^{t_2}\psi_1(t_1)\phi_{j_1}(t_1)dt_1dt_2+\right.
$$

\vspace{2mm}
$$
\left. +
\int\limits_t^T\psi_1(t_2)\phi_{j_2}(t_2)
\int\limits_{t_2}^{T}\psi_2(t_1)\phi_{j_1}(t_1)dt_1\right)dt_2
\phi_{j_1}(t_1)\phi_{j_2}(t_2)
=
$$

\vspace{2mm}
\begin{equation}
\label{334.ye}
=\lim_{n_1,n_2\to\infty}
\sum_{j_1=0}^{n_1}\sum_{j_2=0}^{n_2}\left(C_{j_2j_1}+
C_{j_1j_2}\right)
\phi_{j_1}(t_1)\phi_{j_2}(t_2),\ \ \ (t_1, t_2)\in (t, T)^2.
\end{equation}

\vspace{6mm}

Moreover, the convergence of the series (\ref{334.ye}) is uniform on the
rectangle

\vspace{2mm}
$$
[t+\varepsilon, T-\varepsilon]\times[t+\delta,T-\delta]\ \ \
\hbox{for any} \ \ \ \varepsilon, \delta>0\ \ \ \hbox{(in particular, we can 
choose}\ \varepsilon=\delta).
$$

\vspace{5mm}

In addition, the series (\ref{334.ye}) converges to 
$K'(t_1,t_2)$ at any inner point of the square $[t, T]^2.$
Note that Theorem 4 does not answer the question of convergence
of the series
(\ref{334.ye}) on a boundary of the square
$[t, T]^2.$
In obtaining (\ref{334.ye}) we replaced the order of
integration in the second iterated integral.

Let us substitute $t_1=t_2$ into (\ref{334.ye}).
After that, let us rewrite the limit on 
the right-hand side of (\ref{334.ye}) as two limits.
Let us replace $j_1$ with $j_2$, $j_2$ with $j_1$, 
$n_1$ with $n_2,$ and $n_2$ with $n_1$ in the second limit. Thus, we get

\vspace{-2mm}
\begin{equation}
\label{5656}
\lim\limits_{n_1,n_2\to\infty}
\sum_{j_1=0}^{n_1} \sum_{j_2=0}^{n_2} C_{j_2j_1}\phi_{j_1}(t_1)
\phi_{j_2}(t_1)=\frac{1}{2}\psi_1(t_1)\psi_2(t_1),\ \ \ t_1\in (t, T).
\end{equation}

\vspace{4mm}

According to the above reasoning, the equality
(\ref{5656}) holds uniformly on the interval
$[t+\varepsilon, T-\varepsilon]$ for any $\varepsilon>0.$
Additionally, (\ref{5656}) holds at each interior point
of the interval $[t, T].$

Let us fix $\varepsilon>0$ and 
integrate the equality (\ref{5656}) at the interval $[t+\varepsilon,
T-\varepsilon]$. 
Due to the uniform convergence of the series (\ref{5656}) we can swap  
the series and the integral 

\begin{equation}
\label{5657}
\lim\limits_{n_1,n_2\to\infty}
\sum_{j_1=0}^{n_1} \sum_{j_2=0}^{n_2} C_{j_2j_1}
\int\limits_{t+\varepsilon}^{T-\varepsilon}\phi_{j_1}(t_1)
\phi_{j_2}(t_1)dt_1=\frac{1}{2}\int\limits_{t+\varepsilon}^{T-\varepsilon}
\psi_1(t_1)\psi_2(t_1)dt_1.
\end{equation}

\vspace{4mm}

{\bf Lemma 2.}\ {\it 
Under the conditions of Theorem {\rm 6} the following limit

$$
\lim\limits_{n\to\infty}\sum\limits_{j_1=0}^{n}C_{j_1 j_1}
$$

\vspace{3mm}
\noindent
exists and is finite, where $C_{j_1 j_1}$ 
is defined by {\rm (\ref{ppppa})} for $k=2$ and $j_1=j_2,$
i.e. 

$$
C_{j_1j_1}=\int\limits_t^T \psi_2(t_2)\phi_{j_1}(t_2)
\int\limits_t^{t_2} \psi_1(t_1)\phi_{j_1}(t_1)dt_1 dt_2.
$$
}

\vspace{3mm}

The proof of Lemma~2 will be given further in this section.
Using the equality (\ref{5657})
for $n_1=n_2=n$ and Lemma~2, we get

$$
\frac{1}{2}\int\limits_{t+\varepsilon}^{T-\varepsilon}
\psi_1(t_1)\psi_2(t_1)dt_1=\lim\limits_{n\to\infty}
\sum_{j_1,j_2=0}^{n} C_{j_2j_1}
\int\limits_{t+\varepsilon}^{T-\varepsilon}\phi_{j_1}(t_1)
\phi_{j_2}(t_1)dt_1=
$$

\vspace{1mm}
$$
=\lim\limits_{n\to\infty}
\sum_{j_1,j_2=0}^{n} C_{j_2j_1}\left(
\int\limits_{t}^{T}\phi_{j_1}(t_1)
\phi_{j_2}(t_1)dt_1-\int\limits_{t}^{t+\varepsilon}\phi_{j_1}(t_1)
\phi_{j_2}(t_1)dt_1-\right.
$$

\vspace{1mm}
$$
\left.
-\int\limits_{T-\varepsilon}^{T}\phi_{j_1}(t_1)
\phi_{j_2}(t_1)dt_1\right)=
$$

\vspace{1mm}
$$
=\lim\limits_{n\to\infty}
\sum_{j_1,j_2=0}^{n} C_{j_2j_1}\left(
{\bf 1}_{\{j_1=j_2\}}-
\biggl(\phi_{j_1}(\theta)
\phi_{j_2}(\theta)+\phi_{j_1}(\lambda)\phi_{j_2}(\lambda)\biggr)
\varepsilon\right)=
$$

\vspace{1mm}
\begin{equation}
\label{566}
=
\lim\limits_{n\to\infty}\sum_{j_1=0}^{n}C_{j_1j_1}-
\varepsilon
\lim\limits_{n\to\infty}
\sum_{j_1,j_2=0}^{n}
C_{j_2j_1}\biggl(\phi_{j_1}(\theta)
\phi_{j_2}(\theta)+
\phi_{j_1}(\lambda)
\phi_{j_2}(\lambda)\biggr),
\end{equation}

\vspace{4mm}
\noindent
where $\theta\in [t,t+\varepsilon],$ $\lambda\in [T-\varepsilon,T]$.
In obtaining (\ref{566}) we used the theorem
on the mean value for the Riemann  
integral and orthonormality of the functions
$\phi_{j}(x)$ for $j=0, 1, 2\ldots$

Applying (\ref{566}), we obtain

$$
\varepsilon
\lim\limits_{n\to\infty}
\sum_{j_1,j_2=0}^{n}
C_{j_2j_1}\biggl(\phi_{j_1}(\theta)
\phi_{j_2}(\theta)+
\phi_{j_1}(\lambda)
\phi_{j_2}(\lambda)\biggr)=
$$

\vspace{1mm}
$$
=\lim\limits_{n\to\infty}\sum_{j_1=0}^{n}C_{j_1j_1}-
\lim\limits_{n\to\infty}
\sum_{j_1,j_2=0}^{n} C_{j_2j_1}
\int\limits_{t+\varepsilon}^{T-\varepsilon}\phi_{j_1}(t_1)
\phi_{j_2}(t_1)dt_1,
$$

\vspace{4mm}
\noindent
where the limits

$$
\lim\limits_{n\to\infty}\sum_{j_1=0}^{n}C_{j_1j_1},\ \ \ 
\lim\limits_{n\to\infty}
\sum_{j_1, j_2=0}^{n} C_{j_2j_1}
\int\limits_{t+\varepsilon}^{T-\varepsilon}\phi_{j_1}(t_1)
\phi_{j_2}(t_1)dt_1
$$

\vspace{4mm}
\noindent
exist and are finite (see Lemma~2 and the equality (\ref{5657})). 
This means that the limit 

\vspace{1mm}
$$
\varepsilon
\lim\limits_{n\to\infty}
\sum_{j_1,j_2=0}^{n}
C_{j_2j_1}\biggl(\phi_{j_1}(\theta)
\phi_{j_2}(\theta)+
\phi_{j_1}(\lambda)
\phi_{j_2}(\lambda)\biggr)
$$

\vspace{3mm}
\noindent
also exists and is finite. 

Suppose that 
the following relations

\begin{equation}
\label{rr1-six}
\left\vert
\sum\limits_{j_1,j_2=0}^n
C_{j_2j_1}\phi_{j_2}(T)\phi_{j_1}(T)\right\vert\le K <\infty,\ \ \
\left\vert
\sum\limits_{j_1,j_2=0}^n
C_{j_2j_1}\phi_{j_2}(t)\phi_{j_1}(t)\right\vert\le K<\infty
\end{equation}

\vspace{4mm}
\noindent
are satisfied (the relations (\ref{rr1-six})
will be proved further in this section); constant $K$
does not depend on $n.$

Note that
$$
\left\vert\varepsilon
\lim\limits_{n\to\infty}
\sum_{j_1,j_2=0}^{n}
C_{j_2j_1}\biggl(\phi_{j_1}(\theta)
\phi_{j_2}(\theta)+
\phi_{j_1}(\lambda)
\phi_{j_2}(\lambda)\biggr)\right\vert=
$$

\vspace{1mm}
\begin{equation}
\label{seven-1000}
=\lim\limits_{n\to\infty} \varepsilon  \left\vert
\sum_{j_1,j_2=0}^{n}
C_{j_2j_1}\phi_{j_1}(\theta)
\phi_{j_2}(\theta)+
\sum_{j_1,j_2=0}^{n}
C_{j_2j_1}
\phi_{j_1}(\lambda)
\phi_{j_2}(\lambda)\right\vert.
\end{equation}

\vspace{4mm}

Using (\ref{5656}) ($n_1=n_2=n$) and (\ref{rr1-six}), we obtain

$$
\varepsilon  \left\vert
\sum_{j_1,j_2=0}^{n}
C_{j_2j_1}\phi_{j_1}(\theta)
\phi_{j_2}(\theta)+
\sum_{j_1,j_2=0}^{n}
C_{j_2j_1}
\phi_{j_1}(\lambda)
\phi_{j_2}(\lambda)\right\vert\le
$$

\vspace{1mm}
\begin{equation}
\label{sixth-1}
\le \varepsilon\left(
\left\vert
\sum_{j_1,j_2=0}^{n}
C_{j_2j_1}\phi_{j_1}(\theta)
\phi_{j_2}(\theta)\right\vert+
\left\vert\sum_{j_1,j_2=0}^{n}
C_{j_2j_1}
\phi_{j_1}(\lambda)
\phi_{j_2}(\lambda)\right\vert\right)\le 2\varepsilon K_1\ \to\ 0
\end{equation}

\vspace{4mm}
\noindent
if $\varepsilon\to +0,$ where 
$\theta\in [t,t+\varepsilon],$ $\lambda\in [T-\varepsilon,T]$,
constant $K_1$ is independent on $n.$

Performing the 
passage to the limit $\lim\limits_{\varepsilon\to +0}$
in the equality (\ref{566}) and taking into account (\ref{seven-1000}), (\ref{sixth-1}), we get 

\vspace{-4mm}
\begin{equation}
\label{seven-1}
\frac{1}{2}\int\limits_t^T\psi_1(t_1)\psi_2(t_1)dt_1
=\sum_{j_1=0}^{\infty}
C_{j_1j_1}.
\end{equation}

\vspace{2mm}

Thus, to complete the proof of Theorem 6, it is necessary to prove 
(\ref{rr1-six}) and Lemma~2.
To prove (\ref{rr1-six}) and Lemma~2, as well as for further consideration, we need 
some well known properties of the Legendre polynomials \cite{suet}, \cite{Gob}.

The complete orthonormal system
of Legendre polynomials in the space $L_2([t, T])$
looks as follows

\begin{equation}
\label{4009}
\phi_j(x)=\sqrt{\frac{2j+1}{T-t}}P_j\biggl(\biggl(
x-\frac{T+t}{2}\biggr)\frac{2}{T-t}\biggr),\ \ \ j=0, 1, 2,\ldots,
\end{equation}

\vspace{4mm}
\noindent
where $P_j(x)$ is the Legendre polynomial. 

It is known that the Legendre polynomial $P_j (x)$ is represented, 
for example, as

$$
P_j(x)=\frac{1}{2^j j!} \frac{d^j}{dx^j}\left(x^2-1\right)^j.
$$

\vspace{4mm}

At the boundary points of the interval of orthogonality the Legendre 
polynomials satisfy the following relations

\vspace{-1mm}
$$
P_j(1)=1,\ \ \ P_j(-1)=(-1)^j, 
$$

\vspace{1mm}
$$
P_{j+1}(1)-P_{j-1}(1)=0,\ \ \ P_{j+1}(-1)-P_{j-1}(-1)=0,
$$

\vspace{1mm}
$$
P_{j+1}(1)-P_j(1)=0,\ \ \ P_{j+1}(-1)+P_j(-1)=0,
$$

\vspace{4mm}
\noindent
where $j=0, 1, 2, \ldots$

Relation of the Legendre polynomial $P_j(x)$ with derivatives 
of the Legendre polynomials $P_{j+1}(x)$ and $P_{j-1}(x)$ is expressed 
by the following equality

\vspace{1mm}
\begin{equation}
\label{seven-3}
P_j(x)=\frac{1}{2j+1}\left(P_{j+1}^{'}(x)-
P_{j-1}^{'}(x)\right),\ \ \
j=1, 2,\ldots
\end{equation}

\vspace{5mm}

The recurrent relation has the form

\vspace{1mm}
$$
xP_{j}(x)=\frac{(j+1)P_{j+1}(x)+jP_{j-1}(x)}{2j+1},\ \ \ j=1, 2,\ldots
$$

\vspace{4mm}

Orthogonality of Legendre polynomial $P_j(x)$ to any polynomial
$Q_k (x)$ of lesser degree we write in the following form

$$
\int\limits_{-1}^1 Q_k(x) P_j(x)dx=0,\ \ \ k=0, 1, 2,\ldots,j-1.
$$

\vspace{4mm}

From the property

$$
\int\limits_{-1}^1 P_k(x)P_j(x)dx=\left\{
\begin{matrix}
0\ \ &\hbox{if}\ \
k\ne j\cr\cr
2/(2j+1)\ \ &\hbox{if}\ \ k=j
\end{matrix}
\right.
$$

\vspace{5mm}
\noindent
it follows that the orthonormal on the interval $[-1,1]$ Legendre polynomials 
determined by the relation

\vspace{1mm}
$$
P_j^{*}(x)=\sqrt{\frac{2j+1}{2}}P_j(x),\ \ \ j=0, 1, 2,\ldots
$$

\vspace{5mm}

It is well known that there is an estimate

\begin{equation}
\label{otit987}
|P_j(y)| <\frac{K}{\sqrt{j+1}(1-y^2)^{1/4}},\ \ \ 
y\in (-1, 1),\ \ \ j=1, 2,\ldots,
\end{equation}

\vspace{5mm}
\noindent
where constant $K$ does not depends on $y$ and $j$.

Moreover,

\vspace{-3mm}
\begin{equation}
\label{q1}
|P_j(x)|\le 1,\ \ \ x\in[-1,1],\ \ \ j=0,\ 1,\ldots
\end{equation}

\vspace{5mm}

The Christoffel--Darboux formula has the form

\begin{equation}
\label{ty}
\sum\limits_{j=0}^n (2j+1)P_j(x)P_j(y)=
(n+1)\frac{P_n(x)P_{n+1}(y)-P_{n+1}(x)P_{n}(y)}{y-x}.
\end{equation}

\vspace{4mm}

Let us prove (\ref{rr1-six}).
From (\ref{ty}) for $x=\pm 1$ we obtain

\vspace{-1mm}
\begin{equation}
\label{yy3-six}
\sum\limits_{j=0}^n (2j+1)P_j(y)=
(n+1)\frac{P_{n+1}(y)-P_{n}(y)}{y-1},
\end{equation}

\begin{equation}
\label{yy4-six}
\sum\limits_{j=0}^n (2j+1)(-1)^{j}P_j(y)=
(n+1)(-1)^n\frac{P_{n+1}(y)+P_{n}(y)}{y+1}.
\end{equation}

\vspace{4mm}

From the other hand (see (\ref{seven-3}))

$$
\sum\limits_{j=0}^n (2j+1)P_j(y)=1+\sum\limits_{j=1}^n (2j+1)P_j(y)=
$$

\vspace{1mm}
$$
=1+\sum\limits_{j=1}^n (P_{j+1}^{'}(y)-P_{j-1}^{'}(y))=
1+\biggl(\sum\limits_{j=1}^n (P_{j+1}(y)-P_{j-1}(y))\biggr)'=
$$

\vspace{1mm}
\begin{equation}
\label{yy2-six}
=1+(P_{n+1}(x)+P_n(x)-x-1)'=(P_n(x)+P_{n+1}(x))'
\end{equation}

\vspace{4mm}
\noindent
and
$$
\sum\limits_{j=0}^n (2j+1)(-1)^j P_j(y)=1+\sum\limits_{j=1}^n (-1)^j
(2j+1)P_j(y)=
$$

\vspace{1mm}
$$
=1+\sum\limits_{j=1}^n (-1)^j (P_{j+1}^{'}(y)-P_{j-1}^{'}(y))=
1+\biggl(\sum\limits_{j=1}^n (-1)^j(P_{j+1}(y)-P_{j-1}(y))\biggr)'=
$$

\vspace{1mm}
\begin{equation}
\label{yy1-six}
=1+((-1)^n(P_{n+1}(x)-P_n(x))-x+1)'=(-1)^n(P_{n+1}(x)-P_{n}(x))'.
\end{equation}

\vspace{6mm}

Applying (\ref{yy3-six})--(\ref{yy1-six}), we get

\vspace{-1mm}
\begin{equation}
\label{uu1-six}
(n+1)\frac{P_{n+1}(y)-P_{n}(y)}{y-1}=(P_n(x)+P_{n+1}(x))',
\end{equation}

\vspace{1mm}
\begin{equation}
\label{uu2-six}
(n+1)\frac{P_{n+1}(y)+P_{n}(y)}{y+1}
=(P_{n+1}(x)-P_{n}(x))'.
\end{equation}

\vspace{4mm}

Let us prove the boundedness of the first sum in (\ref{rr1-six}). 
We have

$$
\sum\limits_{j_1,j_2=0}^n
C_{j_2j_1}\phi_{j_2}(T)\phi_{j_1}(T)=
$$

\vspace{1mm}
$$
=\frac{1}{4}
\sum\limits_{j_2=0}^n
\sum\limits_{j_1=0}^n
(2j_2+1)(2j_1+1)\int\limits_{-1}^{1}
\psi_2(h(y))P_{j_2}(y)
\int\limits_{-1}^{y}
\psi_1(h(y_1))P_{j_1}(y_1)dy_1 dy=
$$

\vspace{2mm}
$$
=\frac{1}{4}
\int\limits_{-1}^{1}
\psi_2(h(y))
\sum\limits_{j_2=0}^n
(2j_2+1)P_{j_2}(y)
\int\limits_{-1}^{y}
\psi_1(h(y_1))\sum\limits_{j_1=0}^n
(2j_1+1)P_{j_1}(y_1)dy_1 dy=
$$

\vspace{2mm}
$$
=\frac{1}{4}
\int\limits_{-1}^{1}
\psi_2(h(y))
\left(\int\limits_{-1}^{y}
\psi_1(h(y_1))
d(P_{n+1}(y_1)+P_{n}(y_1))\right)d(P_{n+1}(y)+P_{n}(y))=
$$

\vspace{2mm}
$$
=\frac{1}{4}
\int\limits_{-1}^{1}
\psi_1(h(y))
\left(\int\limits_{-1}^{y}
\psi_1(h(y_1))
d(P_{n+1}(y_1)+P_{n}(y_1))\right)d(P_{n+1}(y)+P_{n}(y))+
$$

\vspace{2mm}
$$
+\frac{1}{4}
\int\limits_{-1}^{1}
\Delta(h(y))
\left(\int\limits_{-1}^{y}
\psi_1(h(y_1))
d(P_{n+1}(y_1)+P_{n}(y_1))\right)d(P_{n+1}(y)+P_{n}(y))=
$$

\vspace{1mm}
$$
=\frac{1}{4}I_1+\frac{1}{4}I_2,
$$

\vspace{4mm}
\noindent
where
\begin{equation}
\label{191-six}
\Delta(h(y))=\psi_2(h(y))-\psi_1(h(y)),\ \ \ h(y)=\frac{T-t}{2}y+\frac{T+t}{2}.
\end{equation}

\vspace{4mm}

Further,
$$
I_1=\frac{1}{2}
\biggl(\int\limits_{-1}^{1}
\psi_1(h(y))
d(P_{n+1}(y)+P_{n}(y))\biggr)^2=
$$

\vspace{1mm}
$$
=\frac{1}{2}\biggl(2\psi_1(T)-\int\limits_{-1}^{1}
(P_{n+1}(y)+P_{n}(y))\psi_1'(h(y))\frac{T-t}{2}dy\biggr)^2<C_1<\infty,
$$

\vspace{4mm}
\noindent
where $\psi_1'$ is a derivative 
of the function 
$\psi_1$ with respect to the variable $y,$
constant $C_1$ does not depend on $n$.

By the Lagrange formula we obtain

$$
\Delta(h(y))=\psi_2\biggl(\frac{1}{2}(T-t)(y-1)+T\biggr)-
\psi_1\biggl(\frac{1}{2}(T-t)(y-1)+T\biggr)=
$$

\vspace{1mm}
$$
=\psi_2(T)-\psi_1(T)+(y-1)\biggl(\psi_2'(\xi_y)-\psi_1'(\theta_y)\biggr)
\frac{1}{2}(T-t)=
$$

\vspace{1mm}
\begin{equation}
\label{hu-six}
=C_1+\alpha_y(y-1),
\end{equation}

\vspace{4mm}
\noindent
where $|\alpha_y|<\infty$ and $C_1=\psi_2(T)-\psi_1(T).$

Let us substitute (\ref{hu-six}) into the integral $I_2$ 

$$
I_2=I_3+I_4,
$$

\vspace{3mm}
\noindent
where
$$
I_3=\int\limits_{-1}^{1}
\alpha_y(y-1)
\left(\int\limits_{-1}^{y}
\psi_1(h(y_1))
d(P_{n+1}(y_1)+P_{n}(y_1))\right)d(P_{n+1}(y)+P_{n}(y)),
$$

\vspace{1mm}
$$
I_4=
C_1\int\limits_{-1}^{1}
\left(\int\limits_{-1}^{y}
\psi_1(h(y_1))
d(P_{n+1}(y_1)+P_{n}(y_1))\right)d(P_{n+1}(y)+P_{n}(y)).
$$

\vspace{4mm}

Integrating by parts and using (\ref{uu1-six}), we obtain

$$
I_3=
\int\limits_{-1}^{1}
\frac{\alpha_y(y-1)(n+1)(P_{n+1}(y)-P_{n}(y))}{y-1}
\biggl(
\psi_1(h(y))(P_{n+1}(y)+P_{n}(y))-\biggr.
$$

\vspace{1mm}
$$
\biggl.-
\int\limits_{-1}^{y}
(P_{n+1}(y_1)+P_{n}(y_1))\psi_1'(h(y_1))\frac{1}{2}(T-t)dy_1\biggr)dy.
$$

\vspace{4mm}

Applying the etimate
(\ref{otit987}) and taking into account the boundedness
of $\alpha_y$ and $\psi_1'(h(y_1))$, we have
that $|I_3|<\infty.$

Using the integration order replacement in $I_4,$ we get

$$
I_4=
C_1\int\limits_{-1}^{1}\psi_1(h(y_1))
\left(\int\limits_{y_1}^{1}
d(P_{n+1}(y)+P_{n}(y))\right)d(P_{n+1}(y_1)+P_{n}(y_1))=
$$

\vspace{1mm}
$$
=C_1\int\limits_{-1}^{1}\psi_1(h(y_1))d(P_{n+1}(y_1)+P_{n}(y_1))
\int\limits_{-1}^{1}
d(P_{n+1}(y)+P_{n}(y))-
$$

\vspace{1mm}
$$
-C_1\int\limits_{-1}^{1}\psi_1(h(y_1))
\left(\int\limits_{-1}^{y_1}
d(P_{n+1}(y)+P_{n}(y))\right)d(P_{n+1}(y_1)+P_{n}(y_1))=
$$

\vspace{1mm}
$$
=I_5-I_6.
$$

\vspace{4mm}

Consider $I_5$
$$
I_5=2C_1\int\limits_{-1}^{1}\psi_1(h(y_1))d(P_{n+1}(y_1)+P_{n}(y_1))=
$$

\vspace{1mm}
$$
=2C_1\left(2\psi_1(T)-
\int\limits_{-1}^{1}
(P_{n+1}(y_1)+P_{n}(y_1))
\psi_1'(h(y_1))\frac{1}{2}(T-t)dy_1\right).
$$

\vspace{4mm} 

Applying the estimate (\ref{q1}) and using the boundedness of
$\psi_1'(h(y_1))$, we obtain
that $|I_5|<\infty.$

Since (see (\ref{hu-six}))
$$
\psi_1(h(y))=
\psi_1\biggl(\frac{1}{2}(T-t)(y-1)+T\biggr)=
$$

\vspace{1mm}
$$
=\psi_1(T)+(y-1)\psi_1'(\theta_y)
\frac{1}{2}(T-t)=C_2+\beta_y(y-1),
$$

\vspace{4mm}
\noindent
where $|\beta_y|<\infty$ and $C_2=\psi_1(T),$ then

$$
I_6=
C_3\int\limits_{-1}^{1}
\left(\int\limits_{-1}^{y_1}
d(P_{n+1}(y)+P_{n}(y))\right)d(P_{n+1}(y_1)+P_{n}(y_1))+
$$

\vspace{1mm}
$$
+C_1\int\limits_{-1}^{1}\beta_{y_1}(y_1-1)
\left(\int\limits_{-1}^{y_1}
d(P_{n+1}(y)+P_{n}(y))\right)d(P_{n+1}(y_1)+P_{n}(y_1))=
$$

\vspace{1mm}
$$
=\frac{C_3}{2}\left(\int\limits_{-1}^{1}
d(P_{n+1}(y)+P_{n}(y))\right)^2+
$$

\vspace{1mm}
$$
+C_1\int\limits_{-1}^{1}\frac{\beta_{y_1}(y_1-1)
(n+1)(P_{n+1}(y_1)-P_{n}(y_1))}{y_1-1}
\left(\int\limits_{-1}^{y_1}
d(P_{n+1}(y)+P_{n}(y))\right)dy_1=
$$

\vspace{1mm}
$$
=2C_3+
C_1\int\limits_{-1}^{1}\beta_{y_1}
(n+1)(P_{n+1}(y_1)-P_{n}(y_1))
(P_{n+1}(y_1)+P_{n}(y_1))dy_1.
$$

\vspace{4mm}

Using the estimate (\ref{otit987}) and 
taking into account the bounedness of $\beta_{y_1}$, we obtain
that $|I_6|<\infty.$
Thus, the boundedness of the first sum in (\ref{rr1-six}) is proved.

Let us prove the boundedness of the second sum in (\ref{rr1-six}). 
We have

$$
\sum\limits_{j_1,j_2=0}^n
C_{j_2j_1}\phi_{j_2}(t)\phi_{j_1}(t)=
$$

\vspace{1mm}
$$
=\frac{1}{4}
\sum\limits_{j_2=0}^n
\sum\limits_{j_1=0}^n
(2j_2+1)(2j_1+1)(-1)^{j_1+j_2}\int\limits_{-1}^{1}
\psi_2(h(y))P_{j_2}(y)
\int\limits_{-1}^{y}
\psi_1(h(y_1))P_{j_1}(y_1)dy_1 dy=
$$

\vspace{2mm}
$$
=\frac{1}{4}
\int\limits_{-1}^{1}
\psi_2(h(y))
\sum\limits_{j_2=0}^n
(2j_2+1)P_{j_2}(y)(-1)^{j_2}
\int\limits_{-1}^{y}
\psi_1(h(y_1))
\sum\limits_{j_1=0}^n
(2j_1+1)P_{j_1}(y_1)(-1)^{j_1}
dy_1 dy=
$$

\vspace{2mm}
$$
=\frac{(-1)^{2n}}{4}
\int\limits_{-1}^{1}
\psi_2(h(y))
\left(\int\limits_{-1}^{y}
\psi_1(h(y_1))
d(P_{n+1}(y_1)-P_{n}(y_1))\right)
d(P_{n+1}(y)-P_{n}(y))=
$$

\vspace{2mm}
$$
=\frac{1}{4}
\int\limits_{-1}^{1}
\psi_1(h(y))
\left(\int\limits_{-1}^{y}
\psi_1(h(y_1))
d(P_{n+1}(y_1)-P_{n}(y_1))\right)d(P_{n+1}(y)-P_{n}(y))+
$$

\vspace{2mm}
$$
+\frac{1}{4}
\int\limits_{-1}^{1}
\Delta(h(y))
\left(\int\limits_{-1}^{y}
\psi_1(h(y_1))
d(P_{n+1}(y_1)-P_{n}(y_1))\right)d(P_{n+1}(y)-P_{n}(y))=
$$

\vspace{1mm}
$$
=\frac{1}{4}J_1+\frac{1}{4}J_2,
$$

\vspace{4mm}
\noindent
where $\Delta(h(y)),$
$h(y)$ are defined by (\ref{191-six}).

Further,
$$
J_1=\frac{1}{2}
\left(\int\limits_{-1}^{1}
\psi_1(h(y))
d(P_{n+1}(y)-P_{n}(y))\right)^2=
$$

\vspace{1mm}
\begin{equation}
\label{hh1-six}
=\frac{1}{2}\left(2(-1)^{n}\psi_1(t)-\int\limits_{-1}^{1}
(P_{n+1}(y)-P_{n}(y))\psi_1'(h(y))\frac{T-t}{2}dy\right)^2<K_1<\infty,
\end{equation}

\vspace{4mm}
\noindent
where $\psi_1'$ is a derivative of the function
$\psi_1$ with respect to the variable $y,$
constant $K_1$ is independent of $n$.

By the Lagrange formula we obtain

$$
\Delta(h(y))=\psi_2\biggl(\frac{1}{2}(T-t)(y+1)+t\biggr)-
\psi_1\biggl(\frac{1}{2}(T-t)(y+1)+t\biggr)=
$$

\vspace{1mm}
$$
=\psi_2(t)-\psi_1(t)+(y+1)\biggl(\psi_2'(\mu_y)-\psi_1'(\rho_y)\biggr)
\frac{1}{2}(T-t)=
$$

\vspace{1mm}
\begin{equation}
\label{hu11-six}
=K_2+\gamma_y(y+1),
\end{equation}

\vspace{4mm}
\noindent
where $|\gamma_y|<\infty$ and $K_2=\psi_2(t)-\psi_1(t).$

Consider $J_2$
$$
J_2=
\int\limits_{-1}^{1}
\Delta(h(y))
d(P_{n+1}(y)-P_{n}(y))
\int\limits_{-1}^{1}
\psi_1(h(y_1))
d(P_{n+1}(y_1)-P_{n}(y_1))-
$$

\vspace{1mm}
$$
-
\int\limits_{-1}^{1}
\Delta(h(y))
\left(\int\limits_{y}^{1}
\psi_1(h(y_1))
d(P_{n+1}(y_1)-P_{n}(y_1))\right)d(P_{n+1}(y)-P_{n}(y))=
$$

\vspace{1mm}
$$
=J_3J_4-J_5.
$$

\vspace{4mm}

The integral $J_4$ was considered earlier (see $J_1$ and (\ref{hh1-six})), i.e.
it has already been shown that
$|J_4|<\infty$. Analogously, we have that $|J_3|<\infty$.

Let us substitute (\ref{hu11-six}) into the integral $J_5$

$$
J_5=J_6+J_7,
$$

\vspace{3mm}
\noindent
where
$$
J_6=\int\limits_{-1}^{1}
\gamma_y(y+1)
\left(\int\limits_{y}^{1}
\psi_1(h(y_1))
d(P_{n+1}(y_1)-P_{n}(y_1))\right)d(P_{n+1}(y)-P_{n}(y)),
$$

\vspace{1mm}
$$
J_7=
K_2\int\limits_{-1}^{1}
\left(\int\limits_{y}^{1}
\psi_1(h(y_1))
d(P_{n+1}(y_1)-P_{n}(y_1))\right)d(P_{n+1}(y)-P_{n}(y)).
$$

\vspace{4mm}

Integrating by parts and using (\ref{uu2-six}), we get

$$
J_6=
\int\limits_{-1}^{1}
\frac{\gamma_y(y+1)(n+1)(P_{n+1}(y)+P_{n}(y))}{y+1}
\biggl(
-\psi_1(h(y))(P_{n+1}(y)-P_{n}(y))-\biggr.
$$

\vspace{1mm}
$$
\biggl.-
\int\limits_{y}^{1}
(P_{n+1}(y_1)-P_{n}(y_1))\psi_1'(h(y_1))\frac{1}{2}(T-t)dy_1\biggr)dy.
$$

\vspace{4mm}

Applying the etimate
(\ref{otit987}) and taking into account the boundedness
of $\gamma_y$ and $\psi_1'(h(y_1))$, we have
that $|J_6|<\infty.$

Using the integration order replacement in $J_7,$ we obtain

$$
J_7=
K_2\int\limits_{-1}^{1}
\psi_1(h(y_1))
\left(\int\limits_{-1}^{y_1}
d(P_{n+1}(y)-P_{n}(y))\right)d(P_{n+1}(y_1)-P_{n}(y_1))=
$$

\vspace{1mm}
$$
=K_2\int\limits_{-1}^{1}
\psi_1(h(y_1))d(P_{n+1}(y_1)-P_{n}(y_1))
\int\limits_{-1}^{1}
d(P_{n+1}(y)-P_{n}(y))-K_2J_8=
$$

\vspace{3mm}
$$
=K_2 J_4 2(-1)^{n} - K_2J_8,
$$

\vspace{4mm}
\noindent
where
$$
J_8=\int\limits_{-1}^{1}
\psi_1(h(y_1))
\left(\int\limits_{y_1}^{1}
d(P_{n+1}(y)-P_{n}(y))\right)d(P_{n+1}(y_1)-P_{n}(y_1)).
$$

\vspace{4mm}

Since (see (\ref{hu11-six}))
$$
\psi_1(h(y))=\psi_1\biggl(\frac{1}{2}(T-t)(y+1)+t\biggr)=
$$

\vspace{1mm}
\begin{equation}
\label{ooo1}
=\psi_1(t)+(y+1)\psi_1'(\rho_y)
\frac{1}{2}(T-t)=K_3+\varepsilon_y(y+1),
\end{equation}

\vspace{4mm}
\noindent
where $|\varepsilon_y|<\infty$ and $K_3=\psi_1(t),$ then

$$
J_8=K_3\int\limits_{-1}^{1}
\left(\int\limits_{y_1}^{1}
d(P_{n+1}(y)-P_{n}(y))\right)d(P_{n+1}(y_1)-P_{n}(y_1))+
$$

\vspace{1mm}
$$
+\int\limits_{-1}^{1}
\varepsilon_y(y+1)
\left(\int\limits_{y_1}^{1}
d(P_{n+1}(y)-P_{n}(y))\right)d(P_{n+1}(y_1)-P_{n}(y_1))=
$$

\vspace{1mm}
$$
=\frac{K_3}{2}\left(\int\limits_{-1}^{1}
d(P_{n+1}(y)-P_{n}(y))\right)^2+
$$

\vspace{1mm}
$$
+\int\limits_{-1}^{1}
\frac{\varepsilon_{y_1}(y_1+1)(n+1)(P_{n+1}(y_1)+P_{n}(y_1))}{y_1+1}
(P_{n}(y_1)-P_{n+1}(y_1))dy=
$$

\vspace{1mm}
\begin{equation}
\label{qqq1}
=2K_3+\int\limits_{-1}^{1}
\varepsilon_{y_1}(n+1)(P_{n+1}(y_1)+P_{n}(y_1))
(P_{n}(y_1)-P_{n+1}(y_1))dy.
\end{equation}

\vspace{2mm}

When obtaining the equality (\ref{qqq1}), we used (\ref{uu2-six}).
Applying the estimate (\ref{otit987}) and 
taking into account the bounedness of $\varepsilon_{y_1}$, we obtain
that $|J_8|<\infty.$
Thus, the boundedness of the second sum in (\ref{rr1-six}) is proved.
The relations (\ref{rr1-six}) are proved.

Let us prove Lemma~2. We will prove
that 

$$
\sum\limits_{j_1=0}^n C_{j_1j_1}
$$

\vspace{3mm}
\noindent
is the Cauchy sequence for 
the cases of Legendre polynomials and trigonometric functions.

Consider the case of Legendre polynomials
and fix $n>m$ $(n, m\in \mathbb{N})$. We have

\vspace{1mm}
$$
\sum\limits_{j_1=m+1}^n
C_{j_1j_1}=
\sum\limits_{j_1=m+1}^n
\int\limits_t^T \psi_2(s)\phi_{j_1}(s)
\int\limits_t^{s} \psi_1(\tau)\phi_{j_1}(\tau)d\tau ds=
$$

\vspace{3mm}
$$
=
\frac{T-t}{4}
\sum\limits_{j_1=m+1}^n
(2j_1+1)\int\limits_{-1}^{1}
\psi_2(h(x))P_{j_1}(x)
\int\limits_{-1}^{x}
\psi_1(h(y))P_{j_1}(y)dy dx=
$$

\vspace{3mm}
$$
=
\frac{T-t}{4}
\sum\limits_{j_1=m+1}^n
\int\limits_{-1}^{1}
\psi_1(h(x))\psi_2(h(x))\left(P_{j_1+1}(x)P_{j_1}(x)
-P_{j_1}(x)P_{j_1-1}(x)\right)dx-
$$

\vspace{3mm}
$$
-\frac{(T-t)^2}{8}\hspace{-2mm}
\sum\limits_{j_1=m+1}^n
\int\limits_{-1}^{1}
\psi_2(h(x))P_{j_1}(x)
\int\limits_{-1}^{x}
\left(P_{j_1+1}(y)-P_{j_1-1}(y)\right)\psi_1'(h(y))dy dx=
$$

\vspace{3mm}
$$
=
\frac{T-t}{4}
\int\limits_{-1}^{1}
\psi_1(h(x))\psi_2(h(x))\sum\limits_{j_1=m+1}^n
\left(P_{j_1+1}(x)P_{j_1}(x)
-P_{j_1}(x)P_{j_1-1}(x)\right)dx-
$$

\vspace{3mm}
$$
-\frac{(T-t)^2}{8}\hspace{-2mm}
\sum\limits_{j_1=m+1}^n
\int\limits_{-1}^{1}
\left(P_{j_1+1}(y)-P_{j_1-1}(y)\right)\psi_1'(h(y))
\int\limits_{y}^{1}
P_{j_1}(x)\psi_2(h(x))dx dy=
$$

\vspace{3mm}
$$
=
\frac{T-t}{4}
\int\limits_{-1}^{1}
\psi_1(h(x))\psi_2(h(x))
\left(P_{n+1}(x)P_{n}(x)
-P_{m+1}(x)P_{m}(x)\right)dx+
$$

\vspace{3mm}
$$
+\frac{(T-t)^2}{8}
\sum\limits_{j_1=m+1}^n
\frac{1}{2j_1+1}
\int\limits_{-1}^{1}
\left(P_{j_1+1}(y)-P_{j_1-1}(y)\right)\psi_1'(h(y))\times
$$

\vspace{3mm}
$$
\times
\Biggl(
\left(P_{j_1+1}(y)-P_{j_1-1}(y)\right)\psi_2(h(y))+\Biggr.
$$

\vspace{3mm}
\begin{equation}
\label{tupo14}
\Biggl.
+
\frac{T-t}{2}
\int\limits_{y}^{1}
\left(P_{j_1+1}(x)-
P_{j_1-1}(x)\right)\psi_2'(h(x))dx\Biggr)dy,
\end{equation}

\vspace{5mm}
\noindent
where $\psi_1'$, $\psi_2'$ are
derivatives of the functions $\psi_1(\tau)$, $\psi_2(\tau)$ with respect 
to the variable
$h(y)$ (see (\ref{ter1})).

Applying the estimate (\ref{otit987}) and tak\-ing into account 
the boundedness of the functions $\psi_1(\tau)$, $\psi_2(\tau)$
and their derivatives, we finally obtain

\vspace{1mm}
$$
\left\vert\sum\limits_{j_1=m+1}^n
C_{j_1j_1}\right\vert\le
C_1\left(\frac{1}{n}+\frac{1}{m}\right)
\int\limits_{-1}^1 \frac{dx}{\left(1-x^2\right)^{1/2}}+
$$

\vspace{3mm}
$$
+C_2 \sum\limits_{j_1=m+1}^n \frac{1}{j_1^2}\left(
\int\limits_{-1}^1 \frac{dy}{\left(1-y^2\right)^{1/2}}+
\int\limits_{-1}^1 \frac{1}{\left(1-y^2\right)^{1/4}}
\int\limits_{y}^1 \frac{dx}{\left(1-x^2\right)^{1/4}}dy\right)\le
$$

\vspace{3mm}
\begin{equation}
\label{tupo15}
\le C_3\left(\frac{1}{n}+\frac{1}{m}+\sum\limits_{j_1=m+1}^n 
\frac{1}{j_1^2}\right) \to 0
\end{equation}

\vspace{4mm}
\noindent
if $n, m\to\infty$ $(n>m),$
where constants $C_1, C_2, C_3$ do not depend on $n$ and $m$.

Consider the trigonometric case. 
Below in this section
we write $\lim\limits_{n,m\to\infty}$ instead of 
$\lim\limits_{\stackrel{n,m\to\infty}{{}_{n>m}}}.$
Fix $n>m$ $(n, m\in {\bf N}).$ 
Denote

\vspace{1mm}
$$
S_{n,m}\stackrel{\sf def}{=}
\sum\limits_{j_1=m+1}^n
C_{j_1j_1}=
\sum\limits_{j_1=m+1}^n
\int\limits_t^T \psi_2(t_2)\phi_{j_1}(t_2)
\int\limits_t^{t_2} \psi_1(t_1)\phi_{j_1}(t_1)dt_1 dt_2.
$$

\vspace{4mm}

By analogy with (\ref{tupo14})
we obtain

\vspace{1mm}
$$
S_{2n,2m}=
\sum\limits_{j_1=2m+1}^{2n}
\int\limits_t^T \psi_2(t_2)\phi_{j_1}(t_2)
\int\limits_t^{t_2} \psi_1(t_1)\phi_{j_1}(t_1)dt_1 dt_2=
$$

\vspace{3mm}
$$
=\frac{2}{T-t}\sum\limits_{j_1=m+1}^n\left(
\int\limits_t^T \psi_2(t_2){\rm sin}\frac{2\pi j_1(t_2-t)}{T-t}
\int\limits_t^{t_2}\psi_1(t_1){\rm sin}\frac{2\pi j_1(t_1-t)}{T-t}dt_1 dt_2+
\right.
$$

\vspace{3mm}
$$
\left.+
\int\limits_t^T \psi_2(t_2){\rm cos}\frac{2\pi j_1(t_2-t)}{T-t}
\int\limits_t^{t_2}\psi_1(t_1){\rm cos}\frac{2\pi j_1(t_1-t)}{T-t}dt_1 dt_2\right)=
$$

\vspace{3mm}
$$
=\frac{T-t}{2\pi^2}\sum\limits_{j_1=m+1}^n\frac{1}{j_1^2}
\left(\psi_1(t)\left(\psi_2(t)-\psi_2(T)+
\int\limits_t^{T}\psi_2'(t_2){\rm cos}\frac{2\pi j_1(t_2-t)}{T-t}dt_2\right)-\right.
$$

\vspace{3mm}
$$
-\int\limits_t^{T}\psi_1'(t_1){\rm cos}\frac{2\pi j_1(t_1-t)}{T-t}\Biggl(
\psi_2(T)-\psi_2(t_1){\rm cos}\frac{2\pi j_1(t_1-t)}{T-t}-\Biggr.
$$

\vspace{2mm}
$$
\Biggl.-\int\limits_{t_1}^{T}\psi_2'(t_2){\rm cos}\frac{2\pi j_1(t_2-t)}{T-t}dt_2\Biggr)dt_1+
$$

\vspace{2mm}
$$
+\int\limits_t^{T}\psi_1'(t_1){\rm sin}\frac{2\pi j_1(t_1-t)}{T-t}\Biggl(
\psi_2(t_1){\rm sin}\frac{2\pi j_1(t_1-t)}{T-t}+\Biggr.
$$

\vspace{2mm}
\begin{equation}
\label{agentaa1000}
\left.\Biggl.
+\int\limits_{t_1}^{T}\psi_2'(t_2){\rm sin}\frac{2\pi j_1(t_2-t)}{T-t}dt_2\Biggr)dt_1\right),
\end{equation}

\vspace{5mm}
\noindent
where $\psi_1'(\tau),$ $\psi_2'(\tau)$ are
derivatives of the functions $\psi_1(\tau),$ $\psi_2(\tau)$ with respect to the variable
$\tau$.

From (\ref{agentaa1000}) we get

\vspace{-2mm}
\begin{equation}
\label{agentaa1001}
\left|S_{2n,2m}\right|\le
C \sum\limits_{j_1=m+1}^n 
\frac{1}{j_1^2} \to 0
\end{equation}

\vspace{4mm}
\noindent
if $n, m\to\infty\ (n>m),$
where constant $C$ does not depend on $n$ and $m$.

Further,

\vspace{-1mm}
$$
S_{2n-1,2m}=S_{2n,2m}-
$$

\begin{equation}
\label{agentaa1002}
-\frac{2}{T-t}
\int\limits_t^T \psi_2(t_2){\rm cos}\frac{2\pi n(t_2-t)}{T-t}
\int\limits_t^{t_2}\psi_1(t_1){\rm cos}\frac{2\pi n(t_1-t)}{T-t}dt_1 dt_2,
\end{equation}

\vspace{3mm}
$$
S_{2n,2m-1}=S_{2n,2m}+
$$

\begin{equation}
\label{agentaa1003}
+\frac{2}{T-t}
\int\limits_t^T \psi_2(t_2){\rm cos}\frac{2\pi m(t_2-t)}{T-t}
\int\limits_t^{t_2}\psi_1(t_1){\rm cos}\frac{2\pi m(t_1-t)}{T-t}dt_1 dt_2,
\end{equation}

\vspace{3mm}
$$
S_{2n-1,2m-1}=S_{2n,2m-1}-
$$

\vspace{3mm}
$$
-\frac{2}{T-t}
\int\limits_t^T \psi_2(t_2){\rm cos}\frac{2\pi n(t_2-t)}{T-t}
\int\limits_t^{t_2}\psi_1(t_1){\rm cos}\frac{2\pi n(t_1-t)}{T-t}dt_1 dt_2=
$$

\vspace{3mm}
$$
=S_{2n,2m}+
\frac{2}{T-t}
\int\limits_t^T \psi_2(t_2){\rm cos}\frac{2\pi m(t_2-t)}{T-t}
\int\limits_t^{t_2}\psi_1(t_1){\rm cos}\frac{2\pi m(t_1-t)}{T-t}dt_1 dt_2-
$$

\vspace{3mm}
\begin{equation}
\label{agentaa1004}
-\frac{2}{T-t}
\int\limits_t^T \psi_2(t_2){\rm cos}\frac{2\pi n(t_2-t)}{T-t}
\int\limits_t^{t_2}\psi_1(t_1){\rm cos}\frac{2\pi n(t_1-t)}{T-t}dt_1 dt_2.
\end{equation}

\vspace{5mm}

Integrating by parts in (\ref{agentaa1002})--(\ref{agentaa1004}), we obtain

\begin{equation}
\label{agentaa1005}
\left|S_{2n-1,2m}\right|\le \left|S_{2n,2m}\right|+\frac{C_1}{n},
\end{equation}

\vspace{2mm}
\begin{equation}
\label{agentaa1006}
\left|S_{2n,2m-1}\right|\le \left|S_{2n,2m}\right|+\frac{C_1}{m},
\end{equation}

\vspace{2mm}
\begin{equation}
\label{agentaa1007}
\left|S_{2n-1,2m-1}\right|\le \left|S_{2n,2m}\right|+C_1\left(\frac{1}{m}+\frac{1}{n}\right),
\end{equation}

\vspace{5mm}
\noindent
where constant $C_1$ does not depend on $n$ and $m$.

The relations (\ref{agentaa1001}), (\ref{agentaa1005})--(\ref{agentaa1007}) imply that

\begin{equation}
\label{agentaa1008}
\lim\limits_{n,m\to\infty}\left|S_{2n,2m}\right|=
\lim\limits_{n,m\to\infty}\left|S_{2n-1,2m}\right|=
\lim\limits_{n,m\to\infty}\left|S_{2n,2m-1}\right|=
\lim\limits_{n,m\to\infty}\left|S_{2n-1,2m-1}\right|=0.
\end{equation}

\vspace{4mm}

From (\ref{agentaa1008}) we get

\vspace{-3mm}
\begin{equation}
\label{agentaa1009}
\lim\limits_{n,m\to\infty}\left|S_{n,m}\right|=0.
\end{equation}

\vspace{5mm}

Lemma~2 is proved. Theorem~6 is proved.

\vspace{5mm}

\section{Proof of the Equality (\ref{5t11}). The Case of 
an Arbitrary Complete Orthonormal System of Functions in the Space $L_2([t, T])$
and $\psi_1(\tau),\psi_2(\tau)\in L_2([t, T])$}

\vspace{5mm}

{\bf Theorem~7.}\ {\it Suppose that 
$\{\phi_j(x)\}_{j=0}^{\infty}$ is an arbitrary complete orthonormal system of 
functions in the space $L_2([t, T])$ and $\psi_1(\tau), \psi_2(\tau) \in L_2([t, T]).$
Then the following equality

\vspace{-1mm}
\begin{equation}
\label{trace1}
\sum_{j=0}^{\infty}\int\limits_t^T
\psi_2(t_2)\phi_j(t_2)
\int\limits_t^{t_2}
\psi_1(t_1)\phi_j(t_1)dt_1 dt_2=
\frac{1}{2}\int\limits_t^T 
\psi_1(\tau) \psi_2(\tau) d\tau
\end{equation}

\vspace{3mm}
\noindent
is fulfilled.}

\vspace{2mm}

{\bf Proof.}\ First consider the case $\psi_1(\tau)\equiv \psi_2(\tau)$ 
or 

\vspace{-1mm}
\begin{equation}
\label{trace2}
\psi_1(\tau)=\psi_2(\tau)\int\limits_t^{\tau} g(\theta)d\theta,
\end{equation}

\vspace{3mm}
\noindent
where $\tau\in [t, T]$ and $\psi_1(\tau), \psi_2(\tau), g(\tau)\in L_2([t, T]).$

First suppose that 
$\{\phi_j(x)\}_{j=0}^{\infty}$ is an arbitrary complete orthonormal system of 
functions in the space $L_2([t, T])$
such that $\phi_j(x)$ for $j<\infty$
is continuous at the interval $[t, T]$ except may be
for the finite number of points 
of the finite discontinuity.
Furthermore, let $\psi_1(\tau)\equiv \psi_2(\tau)$ 
or the equality (\ref{trace2}) is satisfied.
Here we suppose that
$\psi_1(\tau), \psi_2(\tau), g(\tau)$  
are continuous functions at the interval $[t, T].$

Using Fubini's Theorem, Lebesgue's Dominated Convergence Theorem
and Parseval's equality, we have (see (\ref{trace2}))

$$
\sum_{j=0}^{\infty}\int\limits_t^T
\psi_2(t_2)\phi_j(t_2)
\int\limits_t^{t_2}
\psi_1(t_1)\phi_j(t_1)dt_1 dt_2=
$$

\vspace{2mm}
\begin{equation}
\label{july10000}
=\sum_{j=0}^{\infty}\int\limits_t^T
\psi_2(t_2)\phi_j(t_2)
\int\limits_t^{t_2}
\psi_2(t_1)\phi_j(t_1)
\int\limits_t^{t_1} g(\tau)d\tau dt_1 dt_2=
\end{equation}

\vspace{2mm}

$$
=\sum_{j=0}^{\infty}\int\limits_t^T
g(\tau)\int\limits_{\tau}^{T} \psi_2(t_1)\phi_j(t_1)
\int\limits_{t_1}^{T}
\psi_2(t_2)\phi_j(t_2)
dt_2 dt_1 d\tau =
$$

\vspace{2mm}

\begin{equation}
\label{after001}
=\frac{1}{2}\sum_{j=0}^{\infty}\int\limits_t^T
g(\tau)\left(\int\limits_{\tau}^{T} \psi_2(t_1)\phi_j(t_1)
dt_1\right)^2 d\tau =
\end{equation}

\vspace{2mm}

\begin{equation}
\label{after002}
=\frac{1}{2}\int\limits_t^T 
g(\tau)\sum_{j=0}^{\infty} \left(\int\limits_{t}^{T} {\bf 1}_{\{\tau<t_1\}}\psi_2(t_1)\phi_j(t_1)
dt_1\right)^2 d\tau =
\end{equation}

\vspace{2mm}

$$
=
\frac{1}{2}\int\limits_t^T 
g(\tau)\int\limits_{t}^{T} {\bf 1}_{\{\tau<t_1\}}\psi_2^2(t_1)
dt_1 d\tau =
$$

\vspace{2mm}

$$
=
\frac{1}{2}\int\limits_t^T 
g(\tau)\int\limits_{\tau}^{T} \psi_2^2(t_1)
dt_1 d\tau =
$$

\vspace{2mm}

\begin{equation}
\label{after003}
=
\frac{1}{2}\int\limits_t^T 
\psi_2^2(t_1) \int\limits_{t}^{t_1} g(\tau)
d\tau dt_1 =
\end{equation}

\vspace{2mm}

\begin{equation}
\label{trace10}
=
\frac{1}{2}\int\limits_t^T 
\psi_1(t_1) \psi_2(t_1) dt_1,
\end{equation}

\vspace{4mm}
\noindent 
where the transition from (\ref{after001}) to (\ref{after002})
is based on Lebesgue's 
Dominated Convergence Theorem. The integrable majorant exists due to 
Parseval's equality

$$
\left\vert g(\tau)\right\vert
\sum_{j=0}^{q} \left(\int\limits_{\tau}^{T}\psi_2(t_1)\phi_j(t_1)
dt_1\right)^2\le
\left\vert g(\tau)\right\vert
\sum_{j=0}^{\infty} \left(\int\limits_{t}^{T} {\bf 1}_{\{\tau<t_1\}}\psi_2(t_1)\phi_j(t_1)
dt_1\right)^2=
$$

\vspace{2mm}
$$
=\left\vert g(\tau)\right\vert
\int\limits_{t}^{T}\left({\bf 1}_{\{\tau<t_1\}}\right)^2 \psi_2^2(t_1)
dt_1\le
\left\vert g(\tau)\right\vert \left\Vert \psi_2 \right\Vert_{L_2([t,T])}^2
=C\left\vert g(\tau)\right\vert
$$

\vspace{4mm}
\noindent
almost everywhere on $[t, T]$ with respect to 
Lebesgue's measure, where constant $C$ does not depend on $p.$

From the other hand, using Fubini's Theorem 
and the generalized Parseval equality as well as the transition from (\ref{july10000})
to (\ref{after003}), we get

$$
\sum_{j=0}^{\infty}\int\limits_t^T
\psi_1(t_2)\phi_j(t_2)
\int\limits_t^{t_2}
\psi_2(t_1)\phi_j(t_1)dt_1 dt_2=
$$

\vspace{2mm}
$$
=\sum_{j=0}^{\infty}\int\limits_t^T
\psi_2(t_2)\phi_j(t_2)
\int\limits_t^{t_2}
g(\tau)d\tau
\int\limits_t^{t_2}  
\psi_2(t_1)\phi_j(t_1)
dt_1 dt_2=
$$

\vspace{2mm}

$$
=\sum_{j=0}^{\infty}\int\limits_t^T
\psi_2(t_1)\phi_j(t_1)
\int\limits_{t_1}^{T}
\psi_2(t_2)\phi_j(t_2)
\int\limits_t^{t_2}  
g(\tau)d\tau
dt_2 dt_1=
$$

\vspace{2mm}

$$
=\sum_{j=0}^{\infty}\int\limits_t^T
\psi_2(t_1)\phi_j(t_1) dt_1
\int\limits_{t}^{T}
\psi_2(t_2)\phi_j(t_2)
\int\limits_t^{t_2}  
g(\tau)d\tau
dt_2 -
$$

\vspace{2mm}

$$
-\sum_{j=0}^{\infty}\int\limits_t^T
\psi_2(t_1)\phi_j(t_1)
\int\limits_{t}^{t_1}
\psi_2(t_2)\phi_j(t_2)
\int\limits_t^{t_2}  
g(\tau)d\tau
dt_2 dt_1=
$$

\vspace{2mm}

$$
=
\int\limits_t^T 
\psi_2(t_1) \cdot \psi_2(t_1) \int\limits_{t}^{t_1} g(\tau)
d\tau dt_1 - \frac{1}{2}\int\limits_t^T 
\psi_2^2(t_1) \int\limits_{t}^{t_1} g(\tau)
d\tau dt_1 =
$$

\vspace{2mm}

\begin{equation}
\label{trace11}
=
\frac{1}{2}\int\limits_t^T 
\psi_2^2(t_1) \int\limits_{t}^{t_1} g(\tau)
d\tau dt_1 =
\frac{1}{2}\int\limits_t^T 
\psi_1(t_1) \psi_2(t_1) dt_1.
\end{equation}

\vspace{4mm}

In addition, for the case $\psi_1(\tau)\equiv \psi_2(\tau)$, 
using the Parseval equality, we obtain

$$
\sum_{j=0}^{\infty}\int\limits_t^T
\psi_1(t_2)\phi_j(t_2)
\int\limits_t^{t_2}
\psi_1(t_1)\phi_j(t_1) dt_1 dt_2=
$$

\vspace{2mm}

$$
=\frac{1}{2}\sum_{j=0}^{\infty}
\left(\int\limits_{t}^{T} \psi_1(t_1)\phi_j(t_1)
dt_1\right)^2 =
$$

\begin{equation}
\label{trace12}
=
\frac{1}{2}
\int\limits_{t}^{T} \psi_1^2(t_1)dt_1.
\end{equation}

\vspace{3mm}

The equality (\ref{trace1}) is proved for $\psi_1(\tau)\equiv \psi_2(\tau)$ 
or when the equality (\ref{trace2}) is satisfied.

Suppose that

\vspace{-3mm}
$$
\psi_2(\tau)=(\tau-t)^l,\ \ \ g(\tau)=k (\tau-t)^{k-1},
$$

\vspace{4mm}
\noindent
where $l=0,1,2,\ldots,$\ $k=1,2,\ldots $

From (\ref{trace2}) we have

\vspace{-1mm}
$$
\psi_1(\tau)=\psi_2(\tau)\int\limits_t^{\tau} g(\theta)d\theta=
k(\tau-t)^l \int\limits_t^{\tau} (\theta-t)^{k-1}d\theta=(\tau-t)^{l+k}.
$$

\vspace{4mm}

Taking into account (\ref{trace10})--(\ref{trace12}),
we obtain 

$$
\sum_{j=0}^{\infty}\int\limits_t^T 
(t_2-t)^l \phi_j(t_2)
\int\limits_t^{t_2} 
(t_1-t)^{l+k} \phi_j(t_1)dt_1 dt_2=
$$

\vspace{2mm}

$$
=
\sum_{j=0}^{\infty}\int\limits_t^T 
(t_2-t)^{l+k} \phi_j(t_2)
\int\limits_t^{t_2} 
(t_1-t)^{l} \phi_j(t_1)dt_1 dt_2=
$$

\vspace{2mm}

\begin{equation}
\label{dsds1}
=
\frac{1}{2}\int\limits_t^T (\tau-t)^{2l+k} d\tau,
\end{equation}

\vspace{4mm}
\noindent
where $k, l=0,1,2,\ldots $

The equality similar to (\ref{dsds1})
was obtained in \cite{rybakov7000x} using other arguments.
In addition, the formula similar to (\ref{dsds1}) was used 
in \cite{rybakov7000x} to generalize the equality (\ref{trace1}) 
to the case of an arbitrary 
complete ortho\-nor\-mal system of functions in the space $L_2([t, T])$
and $\psi_1(\tau),\psi_2(\tau)$ $\in $ $L_2([t, T]).$
Consider this approach \cite{rybakov7000x} in more detail.

Let us rewrite the equality (\ref{dsds1}) in the following form

\vspace{-1mm}
\begin{equation}
\label{strange902}
\sum_{j=0}^{\infty}\int\limits_t^T 
(t_2-t)^l \phi_j(t_2)
\int\limits_t^{t_2} 
(t_1-t)^{m} \phi_j(t_1)dt_1 dt_2=
\frac{1}{2}\int\limits_t^T (\tau-t)^{l}(\tau-t)^{m} d\tau,
\end{equation}

\vspace{3mm}
\noindent
where $l, m=0,1,2,\ldots $

Since the equality (\ref{strange902}) is valid for monomials 
with respect to $\tau-t$ ($\tau\in [t, T]$), it will obviously
also be valid for Legendre polynomials that form a complete 
orthonormal system of functions in the space $L_2([t, T]$
and finite linear combinations of Legendre polynomials.

Let $\psi_1(\tau), \psi_2(\tau)\in L_2([t, T])$ and 
$\psi_1^{(p)}(\tau), \psi_2^{(q)}(\tau)$
be approximations of the functions $\psi_1(\tau),\psi_2(\tau)$,
respectively, which are partial sums of the corresponding
Fourier--Legendre series. Then we have (see (\ref{strange902}))

\vspace{-3mm}
\begin{equation}
\label{strange903}
\sum_{j=0}^{\infty}\int\limits_t^T 
\psi_2^{(q)}(t_2) \phi_j(t_2)
\int\limits_t^{t_2} 
\psi_1^{(p)}(\tau) \phi_j(t_1)dt_1 dt_2=
\frac{1}{2}\int\limits_t^T \psi_1^{(p)}(\tau) \psi_2^{(q)}(\tau) d\tau,
\end{equation}

\vspace{4mm}
\noindent
where $p, q\in \mathbb{N},$ the series converges absolutly and 
its sum does not depend on a basis system $\left\{\phi_j(x)\right\}_{j=0}^{\infty}$.

Let us fix $q$ in (\ref{strange903}). 
The right-hand side of (\ref{strange903}) for a fixed $q$ defines
(as a scalar product in $L_2([t, T])$) a linear bounded
(and therefore continuous) functional in $L_2([t, T]),$
which is given by the function $\psi_2^{(q)}$.
The left-hand side of the equality (\ref{strange903})
has the same properties.
Let us implement the passage to the limit $\lim\limits_{p\to\infty}$
in (\ref{strange903})

\vspace{-1mm}
\begin{equation}
\label{strange904}
\sum_{j=0}^{\infty}\int\limits_t^T 
\psi_2^{(q)}(t_2) \phi_j(t_2)
\int\limits_t^{t_2} 
\psi_1(\tau) \phi_j(t_1)dt_1 dt_2=
\frac{1}{2}\int\limits_t^T \psi_1(\tau) \psi_2^{(q)}(\tau) d\tau,
\end{equation}

\vspace{3mm}
\noindent
where $q\in \mathbb{N}.$ The equality (\ref{strange904})
defines a linear bounded functional in $L_2([t, T])$ given by the function 
$\psi_1$. Let us implement the passage to the limit $\lim\limits_{q\to\infty}$
in (\ref{strange904})

\vspace{-1mm}
$$
\sum_{j=0}^{\infty}\int\limits_t^T 
\psi_2(t_2) \phi_j(t_2)
\int\limits_t^{t_2} 
\psi_1(\tau) \phi_j(t_1)dt_1 dt_2=
\frac{1}{2}\int\limits_t^T \psi_1(\tau) \psi_2(\tau) d\tau,
$$

\vspace{3mm}
\noindent
where $\left\{\phi_j(x)\right\}_{j=0}^{\infty}$
is an arbitrary complete orthonormal system of 
functions in the space $L_2([t,T])$ and
$\psi_1(\tau),\psi_2(\tau)\in $ $L_2([t, T]).$

Thus we have the following theorem.

\vspace{2mm}

{\bf Theorem 8}\ \cite{20xxxxx}.\ {\it Suppose that 
$\{\phi_j(x)\}_{j=0}^{\infty}$ is an arbitrary complete orthonormal system of 
functions in the space $L_2([t, T])$ and
$\psi_1(\tau), \psi_2(\tau)\in L_2([t, T]).$ 
Then 
for the iterated Stratonovich stochastic integral

$$
J^{*}[\psi^{(2)}]_{T,t}={\int\limits_t^{*}}^T\psi_2(t_2)
{\int\limits_t^{*}}^{t_2}\psi_1(t_1)d{\bf f}_{t_1}^{(i_1)}
d{\bf f}_{t_2}^{(i_2)}\ \ \ (i_1, i_2=1,\ldots,m)
$$

\vspace{3mm}
\noindent
the following expansion  

\begin{equation}
\label{trace20}
J^{*}[\psi^{(2)}]_{T,t}=\hbox{\vtop{\offinterlineskip\halign{
\hfil#\hfil\cr
{\rm l.i.m.}\cr
$\stackrel{}{{}_{p_1,p_2\to \infty}}$\cr
}} }\sum_{j_1=0}^{p_1}\sum_{j_2=0}^{p_2}
C_{j_2j_1}\zeta_{j_1}^{(i_1)}\zeta_{j_2}^{(i_2)}
\end{equation}

\vspace{4mm}
\noindent
that converges in the mean-square
sence is valid, where the notations are the same as in Theorem {\rm 3.}
}

In Theorem~8, we use a different definition of the Stratonovich stochastic
integral compared to the definition from \cite{KlPl2}. Namely, we use
the definition from \cite{20xxxxx} (Sect.~2.18).

\vspace{5mm}

\section{Some Recent Results on Expansions of Iterated 
Stratonovich Stochastic Integrals of Multiplicities 3 to 6}

\vspace{5mm}

Recently, a new approach to the expansion and mean-square 
approximation of iterated Stratonovich stochastic integrals has been obtained
\cite{20xxxxx} (Sect.~2.10--2.16), \cite{21} (Sect.~13--19), 
\cite{25} (Sect.~5--11), \cite{arxiv-4} (Sect.~7--13), \cite{new-art-1xxy}.
Let us formulate four theorems that were obtained using this approach.

\vspace{2mm}

{\bf Theorem 9}\ \cite{20xxxxx}, \cite{21}, 
\cite{25}, \cite{arxiv-4}.\
{\it Suppose 
that $\{\phi_j(x)\}_{j=0}^{\infty}$ is a complete orthonormal system of 
Legendre polynomials or trigonometric functions in the space $L_2([t, T]).$
Furthermore, let $\psi_1(\tau), \psi_2(\tau),$ $\psi_3(\tau)$ are continuously dif\-ferentiable 
nonrandom functions on $[t, T].$ 
Then for the 
iterated Stra\-to\-no\-vich stochastic integral of third multiplicity

$$
J^{*}[\psi^{(3)}]_{T,t}={\int\limits_t^{*}}^T\psi_3(t_3)
{\int\limits_t^{*}}^{t_3}\psi_2(t_2)
{\int\limits_t^{*}}^{t_2}\psi_1(t_1)
d{\bf w}_{t_1}^{(i_1)}
d{\bf w}_{t_2}^{(i_2)}d{\bf w}_{t_3}^{(i_3)}\ \ \ (i_1,i_2,i_3=0,1,\ldots,m)
$$

\vspace{4mm}
\noindent
the following 
relations

\vspace{-1mm}
\begin{equation}
\label{fin1}
J^{*}[\psi^{(3)}]_{T,t}
=\hbox{\vtop{\offinterlineskip\halign{
\hfil#\hfil\cr
{\rm l.i.m.}\cr
$\stackrel{}{{}_{p\to \infty}}$\cr
}} }
\sum\limits_{j_1, j_2, j_3=0}^{p}
C_{j_3 j_2 j_1}\zeta_{j_1}^{(i_1)}\zeta_{j_2}^{(i_2)}\zeta_{j_3}^{(i_3)},
\end{equation}

\vspace{3mm}
\begin{equation}
\label{fin2}
{\sf M}\left\{\left(
J^{*}[\psi^{(3)}]_{T,t}-
\sum\limits_{j_1, j_2, j_3=0}^{p}
C_{j_3 j_2 j_1}\zeta_{j_1}^{(i_1)}\zeta_{j_2}^{(i_2)}\zeta_{j_3}^{(i_3)}\right)^2\right\}
\le \frac{C}{p}
\end{equation}

\vspace{5mm}
\noindent
are fulfilled, where $i_1, i_2, i_3=0,1,\ldots,m$ in {\rm (\ref{fin1})} and 
$i_1, i_2, i_3=1,\ldots,m$ in {\rm (\ref{fin2})},
constant $C$ is independent of $p,$

$$
C_{j_3 j_2 j_1}=\int\limits_t^T\psi_3(t_3)\phi_{j_3}(t_3)
\int\limits_t^{t_3}\psi_2(t_2)\phi_{j_2}(t_2)
\int\limits_t^{t_2}\psi_1(t_1)\phi_{j_1}(t_1)dt_1dt_2dt_3
$$

\vspace{4mm}
\noindent
and
$$
\zeta_{j}^{(i)}=
\int\limits_t^T \phi_{j}(\tau) d{\bf f}_{\tau}^{(i)}
$$ 

\vspace{2mm}
\noindent
are independent standard Gaussian random variables for various 
$i$ or $j$ {\rm (}in the case when $i\ne 0${\rm );} 
another notations are the same as in Theorems~{\rm 1, 2}.}

\vspace{2mm}

{\bf Theorem 10}\ \cite{20xxxxx}, \cite{21}, 
\cite{25}, \cite{arxiv-4}.\ {\it Let
$\{\phi_j(x)\}_{j=0}^{\infty}$ be a complete orthonormal system of 
Legendre polynomials or trigonometric functions in the space $L_2([t, T]).$
Furthermore, let $\psi_1(\tau), \ldots,$ $\psi_4(\tau)$ be continuously dif\-ferentiable 
nonrandom functions on $[t, T].$ 
Then for the 
iterated Stra\-to\-no\-vich stochastic integral of fourth multiplicity

\begin{equation}
\label{fin0}
J^{*}[\psi^{(4)}]_{T,t}={\int\limits_t^{*}}^T\psi_4(t_4)
{\int\limits_t^{*}}^{t_4}\psi_3(t_3)
{\int\limits_t^{*}}^{t_3}\psi_2(t_2)
{\int\limits_t^{*}}^{t_2}\psi_1(t_1)
d{\bf w}_{t_1}^{(i_1)}
d{\bf w}_{t_2}^{(i_2)}d{\bf w}_{t_3}^{(i_3)}d{\bf w}_{t_4}^{(i_4)}
\end{equation}

\vspace{4mm}
\noindent
the following 
relations

\begin{equation}
\label{fin3}
J^{*}[\psi^{(4)}]_{T,t}
=\hbox{\vtop{\offinterlineskip\halign{
\hfil#\hfil\cr
{\rm l.i.m.}\cr
$\stackrel{}{{}_{p\to \infty}}$\cr
}} }
\sum\limits_{j_1, j_2, j_3,j_4=0}^{p}
C_{j_4j_3 j_2 j_1}\zeta_{j_1}^{(i_1)}\zeta_{j_2}^{(i_2)}\zeta_{j_3}^{(i_3)}\zeta_{j_4}^{(i_4)},
\end{equation}

\vspace{3mm}

\begin{equation}
\label{fin4}
{\sf M}\left\{\left(
J^{*}[\psi^{(4)}]_{T,t}-
\sum\limits_{j_1, j_2, j_3, j_4=0}^{p}
C_{j_4 j_3 j_2 j_1}\zeta_{j_1}^{(i_1)}\zeta_{j_2}^{(i_2)}\zeta_{j_3}^{(i_3)}
\zeta_{j_4}^{(i_4)}
\right)^2\right\}
\le \frac{C}{p^{1-\varepsilon}}
\end{equation}

\vspace{5mm}
\noindent
are fulfilled, where $i_1, \ldots , i_4=0,1,\ldots,m$ in {\rm (\ref{fin0}),} {\rm (\ref{fin3})} 
and $i_1, \ldots, i_4=1,\ldots,m$ in {\rm (\ref{fin4}),}
constant $C$ does not depend on $p,$
$\varepsilon$ is an arbitrary
small positive real number 
for the case of complete orthonormal system of 
Legendre polynomials in the space $L_2([t, T])$
and $\varepsilon=0$ for the case of
complete orthonormal system of 
trigonometric functions in the space $L_2([t, T]),$

$$
C_{j_4 j_3 j_2 j_1}=
$$

$$
=
\int\limits_t^T\psi_4(t_4)\phi_{j_4}(t_4)
\int\limits_t^{t_4}\psi_3(t_3)\phi_{j_3}(t_3)
\int\limits_t^{t_3}\psi_2(t_2)\phi_{j_2}(t_2)
\int\limits_t^{t_2}\psi_1(t_1)\phi_{j_1}(t_1)dt_1dt_2dt_3dt_4;
$$

\vspace{4mm}
\noindent
another notations are the same as in Theorem~{\rm 9}.}

\vspace{2mm}

{\bf Theorem 11}\ \cite{20xxxxx}, \cite{21}, 
\cite{25}, \cite{arxiv-4}.\
{\it Assume 
that $\{\phi_j(x)\}_{j=0}^{\infty}$ is a complete orthonormal system of 
Legendre polynomials or trigonometric functions in the space $L_2([t, T])$
and $\psi_1(\tau), \ldots,$ $\psi_5(\tau)$ are continuously dif\-ferentiable 
nonrandom functions on $[t, T].$ 
Then for the 
iterated Stra\-to\-no\-vich stochastic integral of fifth multiplicity

\begin{equation}
\label{fin7}
J^{*}[\psi^{(5)}]_{T,t}={\int\limits_t^{*}}^T\psi_5(t_5)
\ldots
{\int\limits_t^{*}}^{t_2}\psi_1(t_1)
d{\bf w}_{t_1}^{(i_1)}
\ldots d{\bf w}_{t_5}^{(i_5)}
\end{equation}

\vspace{4mm}
\noindent
the following 
relations

\begin{equation}
\label{fin8}
J^{*}[\psi^{(5)}]_{T,t}
=\hbox{\vtop{\offinterlineskip\halign{
\hfil#\hfil\cr
{\rm l.i.m.}\cr
$\stackrel{}{{}_{p\to \infty}}$\cr
}} }
\sum\limits_{j_1,\ldots,j_5=0}^{p}
C_{j_5 \ldots j_1}\zeta_{j_1}^{(i_1)}\ldots \zeta_{j_5}^{(i_5)},
\end{equation}

\vspace{3mm}

\begin{equation}
\label{fin9}
{\sf M}\left\{\left(
J^{*}[\psi^{(5)}]_{T,t}-
\sum\limits_{j_1, \ldots, j_5=0}^{p}
C_{j_5 \ldots j_1}\zeta_{j_1}^{(i_1)}\ldots
\zeta_{j_5}^{(i_5)}
\right)^2\right\}
\le \frac{C}{p^{1-\varepsilon}}
\end{equation}

\vspace{5mm}
\noindent
are fulfilled, where $i_1, \ldots , i_5=0,1,\ldots,m$ in {\rm (\ref{fin7}),} {\rm (\ref{fin8})} 
and $i_1, \ldots, i_5=1,\ldots,m$ in {\rm (\ref{fin9}),}
constant $C$ is independent of $p,$
$\varepsilon$ is an arbitrary
small positive real number 
for the case of complete orthonormal system of 
Legendre polynomials in the space $L_2([t, T])$
and $\varepsilon=0$ for the case of
complete orthonormal system of 
trigonometric functions in the space $L_2([t, T]),$

$$
C_{j_5 \ldots j_1}=
\int\limits_t^T\psi_5(t_5)\phi_{j_5}(t_5)\ldots
\int\limits_t^{t_2}\psi_1(t_1)\phi_{j_1}(t_1)dt_1\ldots dt_5;
$$

\vspace{3mm}
\noindent
another notations are the same as in Theorems~{\rm 9, 10}.}

\vspace{2mm}

{\bf Theorem 12}\ \cite{20xxxxx}, \cite{21}, 
\cite{25}, \cite{arxiv-4}.\
{\it Suppose that 
$\{\phi_j(x)\}_{j=0}^{\infty}$ is a complete orthonormal system of 
Legendre polynomials or trigonometric functions in the space $L_2([t, T]).$
Then, for the 
iterated Stratonovich stochastic integral of sixth multiplicity

\begin{equation}
\label{after10001qu1}
J_{T,t}^{*(i_1\ldots i_6)}={\int\limits_t^{*}}^T
\ldots
{\int\limits_t^{*}}^{t_2}
d{\bf w}_{t_1}^{(i_1)}
\ldots d{\bf w}_{t_6}^{(i_6)}
\end{equation}

\vspace{3mm}
\noindent
the following 
expansion 

$$
J_{T,t}^{*(i_1\ldots i_6)}
=\hbox{\vtop{\offinterlineskip\halign{
\hfil#\hfil\cr
{\rm l.i.m.}\cr
$\stackrel{}{{}_{p\to \infty}}$\cr
}} }
\sum\limits_{j_1, \ldots, j_6=0}^{p}
C_{j_6 \ldots j_1}\zeta_{j_1}^{(i_1)}\ldots
\zeta_{j_6}^{(i_6)}
$$

\vspace{3mm}
\noindent
that converges in the mean-square sense is valid, where
$i_1, \ldots, i_6=0, 1,\ldots,m,$

$$
C_{j_6 \ldots j_1}=
\int\limits_t^T\phi_{j_6}(t_6)\ldots
\int\limits_t^{t_2}\phi_{j_1}(t_1)dt_1\ldots dt_6;
$$

\vspace{3mm}
\noindent
another notations are the same as in Theorems~{\rm 9--11}.}

The results of Theorems~9--12 were developed in \cite{20xxxxx} (Chapter~2), 
\cite{21}, \cite{25}, \cite{arxiv-4}.
In particular, analogues of Theorem~12 for iterated Stratonovich stochastic
integrals of multiplicities 7 and 8 were obtained in \cite{20xxxxx} (Sect.~2.36, 2.37).
In addition, the variants of Thorems 9--12 
were obtained
for the case when $\{\phi_j(x)\}_{j=0}^{\infty}$ is an arbitrary complete orthonormal system
of functions in $L_2([t, T])$ \cite{20xxxxx} (Sect.~2.1.4, 2.23, 2.24, 2.31--2.34),
\cite{21}, \cite{25}, \cite{arxiv-4}.

\vspace{5mm}

\section{Expansions of Iterated Stratonovich Stochastic Integrals
of First and Second Multiplicity Based on Multiple Fourier--Legendre
Series}

\vspace{5mm}

We will use the following notations for iterated 
Stratonovich stochastic integrals of first and second multiplicities

\vspace{-2mm}
\begin{equation}
\label{ll1}
I_{(l_1)T,t}^{*(i_1)}
={\int\limits_t^{*}}^T(t-t_1)^{l_1} d{\bf f}_{t_1}^{(i_1)},
\end{equation}

\begin{equation}
\label{ll2}
I_{(l_1 l_2)T,t}^{*(i_1i_2)}
={\int\limits_t^{*}}^T(t-t_2)^{l_2}{\int\limits_t^{*}}^{t_{2}}
(t-t_1)^{l_1} d{\bf f}_{t_1}^{(i_1)}d{\bf f}_{t_2}^{(i_2)},
\end{equation}

\vspace{3mm}
\noindent
where $l_1, l_2=0, 1,\ldots;$ $i_1,\ldots,i_k=1,\ldots,m.$

Note that together with the iterated Stratonovich stochastic integrals
of higher multiplicities than the second, the stochastic integrals
(\ref{ll1}) and (\ref{ll2}) are included in the so-called
unified Taylor--Stratonovich expansion \cite{32} (also see \cite{20xxxxx}-\cite{20eee}). 
This expansion can be used for construction of high-order
strong 
numerical methods for Ito SDEs
(definition 
of a strong numerical method see, for example, in \cite{KlPl2}).

Consider the expansions of some stochastic integrals 
(\ref{ll1}) and (\ref{ll2})
obtained by using Theorems 6, 8 (see (\ref{jes}) or (\ref{trace20}))

$$
I_{(0)T,t}^{*(i_1)}=\sqrt{T-t}\zeta_0^{(i_1)},
$$

\vspace{1mm}

\begin{equation}
\label{4002}
I_{(1)T,t}^{*(i_1)}=-\frac{(T-t)^{3/2}}{2}\left(\zeta_0^{(i_1)}+
\frac{1}{\sqrt{3}}\zeta_1^{(i_1)}\right),
\end{equation}

\vspace{2mm}

\begin{equation}
\label{4003}
I_{(2)T,t}^{*(i_1)}=\frac{(T-t)^{5/2}}{3}\left(\zeta_0^{(i_1)}+
\frac{\sqrt{3}}{2}\zeta_1^{(i_1)}+
\frac{1}{2\sqrt{5}}\zeta_2^{(i_1)}\right),
\end{equation}

\vspace{2mm}

\begin{equation}
\label{4004}
I_{(00)T,t}^{*(i_1 i_2)}=
\frac{T-t}{2}\left(\zeta_0^{(i_1)}\zeta_0^{(i_2)}+\sum_{i=1}^{\infty}
\frac{1}{\sqrt{4i^2-1}}\left(
\zeta_{i-1}^{(i_1)}\zeta_{i}^{(i_2)}-
\zeta_i^{(i_1)}\zeta_{i-1}^{(i_2)}\right)\right),
\end{equation}

\vspace{6mm}

$$
I_{(01)T,t}^{*(i_1 i_2)}=-\frac{T-t}{2}I_{(00)T,t}^{*(i_1 i_2)}
-\frac{(T-t)^2}{4}\left(
\frac{\zeta_0^{(i_1)}\zeta_1^{(i_2)}}{\sqrt{3}}+\right.
$$

\vspace{2mm}
$$
+\left.\sum_{i=0}^{\infty}\left(
\frac{(i+2)\zeta_i^{(i_1)}\zeta_{i+2}^{(i_2)}
-(i+1)\zeta_{i+2}^{(i_1)}\zeta_{i}^{(i_2)}}
{\sqrt{(2i+1)(2i+5)}(2i+3)}-
\frac{\zeta_i^{(i_1)}\zeta_{i}^{(i_2)}}{(2i-1)(2i+3)}\right)\right),
$$

\vspace{8mm}

$$
I_{(10)T,t}^{*(i_1 i_2)}=-\frac{T-t}{2}I_{(00)T,t}^{*(i_1 i_2)}
-\frac{(T-t)^2}{4}\left(
\frac{\zeta_0^{(i_2)}\zeta_1^{(i_1)}}{\sqrt{3}}+\right.
$$

\vspace{2mm}
\begin{equation}
\label{4006}
+\left.\sum_{i=0}^{\infty}\left(
\frac{(i+1)\zeta_{i+2}^{(i_2)}\zeta_{i}^{(i_1)}
-(i+2)\zeta_{i}^{(i_2)}\zeta_{i+2}^{(i_1)}}
{\sqrt{(2i+1)(2i+5)}(2i+3)}+
\frac{\zeta_i^{(i_1)}\zeta_{i}^{(i_2)}}{(2i-1)(2i+3)}\right)\right),
\end{equation}

\vspace{8mm}

$$
I_{(02)T,t}^{*(i_1 i_2)}=-\frac{(T-t)^2}{4}I_{(00)T,t}^{*(i_1 i_2)}
-(T-t)I_{(01)T,t}^{*(i_1 i_2)}+
\frac{(T-t)^3}{8}\left(
\frac{2\zeta_2^{(i_2)}\zeta_0^{(i_1)}}{3\sqrt{5}}+\right.
$$

\vspace{2mm}
$$
+\frac{1}{3}\zeta_0^{(i_1)}\zeta_0^{(i_2)}+
\sum_{i=0}^{\infty}\left(
\frac{(i+2)(i+3)\zeta_{i+3}^{(i_2)}\zeta_{i}^{(i_1)}
-(i+1)(i+2)\zeta_{i}^{(i_2)}\zeta_{i+3}^{(i_1)}}
{\sqrt{(2i+1)(2i+7)}(2i+3)(2i+5)}+
\right.
$$

\vspace{2mm}
$$
\left.\left.+
\frac{(i^2+i-3)\zeta_{i+1}^{(i_2)}\zeta_{i}^{(i_1)}
-(i^2+3i-1)\zeta_{i}^{(i_2)}\zeta_{i+1}^{(i_1)}}
{\sqrt{(2i+1)(2i+3)}(2i-1)(2i+5)}\right)\right),
$$

\vspace{8mm}

$$
I_{(20)T,t}^{*(i_1 i_2)}=-\frac{(T-t)^2}{4}I_{(00)T,t}^{*(i_1 i_2)}
-(T-t)I_{(10)T,t}^{*(i_1 i_2)}+
\frac{(T-t)^3}{8}\left(
\frac{2\zeta_0^{(i_2)}\zeta_2^{(i_1)}}{3\sqrt{5}}+\right.
$$

\vspace{2mm}
$$
+\frac{1}{3}\zeta_0^{(i_1)}\zeta_0^{(i_2)}+
\sum_{i=0}^{\infty}\left(
\frac{(i+1)(i+2)\zeta_{i+3}^{(i_2)}\zeta_{i}^{(i_1)}
-(i+2)(i+3)\zeta_{i}^{(i_2)}\zeta_{i+3}^{(i_1)}}
{\sqrt{(2i+1)(2i+7)}(2i+3)(2i+5)}+
\right.
$$

\vspace{2mm}
$$
\left.\left.+
\frac{(i^2+3i-1)\zeta_{i+1}^{(i_2)}\zeta_{i}^{(i_1)}
-(i^2+i-3)\zeta_{i}^{(i_2)}\zeta_{i+1}^{(i_1)}}
{\sqrt{(2i+1)(2i+3)}(2i-1)(2i+5)}\right)\right),
$$

\vspace{10mm}

$$
I_{(11)T,t}^{*(i_1 i_2)}=-\frac{(T-t)^2}{4}I_{(00)T,t}^{*(i_1 i_2)}
-\frac{(T-t)}{2}\left(I_{(10)T,t}^{*(i_1 i_2)}+
I_{(01)T,t}^{*(i_1 i_2)}\right)+
$$

\vspace{2mm}
$$
+
\frac{(T-t)^3}{8}\left(
\frac{1}{3}\zeta_1^{(i_1)}\zeta_1^{(i_2)}
+
\sum_{i=0}^{\infty}\left(
\frac{(i+1)(i+3)\left(\zeta_{i+3}^{(i_2)}\zeta_{i}^{(i_1)}
-\zeta_{i}^{(i_2)}\zeta_{i+3}^{(i_1)}\right)}
{\sqrt{(2i+1)(2i+7)}(2i+3)(2i+5)}+
\right.\right.
$$

\vspace{2mm}
$$
\left.\left.+
\frac{(i+1)^2\left(\zeta_{i+1}^{(i_2)}\zeta_{i}^{(i_1)}
-\zeta_{i}^{(i_2)}\zeta_{i+1}^{(i_1)}\right)}
{\sqrt{(2i+1)(2i+3)}(2i-1)(2i+5)}\right)\right),
$$

\vspace{7mm}

$$
I_{(3)T,t}^{*(i_1)}=-\frac{(T-t)^{7/2}}{4}\left(\zeta_0^{(i_1)}+
\frac{3\sqrt{3}}{5}\zeta_1^{(i_1)}+
\frac{1}{\sqrt{5}}\zeta_2^{(i_1)}+
\frac{1}{5\sqrt{7}}\zeta_3^{(i_1)}\right),
$$

\vspace{6mm}
\noindent
where 
\begin{equation}
\label{zq1}
\zeta_{j}^{(i)}=
\int\limits_t^T \phi_{j}(\tau) d{\bf f}_{\tau}^{(i)}\ \ \ (i=1,\ldots,m)
\end{equation} 

\vspace{3mm}
\noindent
are independent standard Gaussian random variables
for 
various
$i$ or $j$.

Note the simplicity of the formulas (\ref{4002}), (\ref{4003}).
For comparison, we present analogs of the formulas 
(\ref{4002}), (\ref{4003})
obtained in \cite{KPW} (also see \cite{KlPl2}) using the method proposed
in \cite{Mi2}

\vspace{2mm}
\begin{equation}
\label{42}
I_{(1)T,t}^{(i_1)q}=-\frac{{(T-t)}^{3/2}}{2}
\left(\zeta_0^{(i_1)}-\frac{\sqrt{2}}{\pi}\left(\sum_{r=1}^{q}
\frac{1}{r}
\zeta_{2r-1}^{(i_1)}+\sqrt{\alpha_q}\xi_q^{(i_1)}\right)
\right),
\end{equation}

\vspace{5mm}
\begin{equation}
\label{46}
I_{(2)T,t}^{(i_1)q}=
(T-t)^{5/2}\left(
\frac{1}{3}\zeta_0^{(i_1)}+\frac{1}{\sqrt{2}\pi^2}
\left(\sum_{r=1}^{q}\frac{1}{r^2}\zeta_{2r}^{(i_1)}+
\sqrt{\beta_q}
\mu_q^{(i_1)}\right)-\right.
$$
$$
\left.-\frac{1}{\sqrt{2}\pi}\left(\sum_{r=1}^q
\frac{1}{r}\zeta_{2r-1}^{(i_1)}+
\sqrt{\alpha_q}\xi_q^{(i_1)}\right)\right),
\end{equation}

\vspace{5mm}
\noindent
where $\zeta_j^{(i)}$ is defined by the formula (\ref{zq1}),
$\phi_j(s)$ is a complete orthonormal system of trigonometric
functions in the space
$L_2([t, T]),$ and
$\zeta_0^{(i)},$ $\zeta_{2r}^{(i)},$
$\zeta_{2r-1}^{(i)},$ $\xi_q^{(i)},$ $\mu_q^{(i)}$ $(r=1,\ldots,q,$\ \ 
$i=1,\ldots,m)$ are independent 
standard Gaussian random variables, $i_1=1,\ldots,m,$ 

\vspace{2mm}
$$
\xi_q^{(i)}=\frac{1}{\sqrt{\alpha_q}}\sum_{r=q+1}^{\infty}
\frac{1}{r}\zeta_{2r-1}^{(i)},\ \ \ \
\alpha_q=\frac{\pi^2}{6}-\sum_{r=1}^q\frac{1}{r^2},
$$

\vspace{2mm}
$$
\mu_q^{(i)}=\frac{1}{\sqrt{\beta_q}}\sum_{r=q+1}^{\infty}
\frac{1}{r^2}~\zeta_{2r}^{(i)},\ \ \ \
\beta_q=\frac{\pi^4}{90}-\sum_{r=1}^q\frac{1}{r^4}.
$$

\vspace{6mm}

Another example of obvious advantage of the Legendre polynomials 
over the trigonometric func\-ti\-ons (in the framework of the considered 
problem) is the truncated expansion 
of the iterated Stra\-to\-no\-vich stochastic integral
$I_{(10)T, t}^{*(i_1 i_2)}$ obtained by Theorems 6, 8 in which
instead of the double Fourier--Le\-gen\-dre series is taken 
the double trigonometric Fourier series

$$
I_{(10)T,t}^{*(i_1 i_2)q}=-(T-t)^{2}\Biggl(\frac{1}{6}
\zeta_{0}^{(i_1)}\zeta_{0}^{(i_2)}-\frac{1}{2\sqrt{2}\pi}
\sqrt{\alpha_q}\xi_q^{(i_2)}\zeta_0^{(i_1)}+\Biggr.
$$

\vspace{2mm}
$$
+\frac{1}{2\sqrt{2}\pi^2}\sqrt{\beta_q}\Biggl(
\mu_q^{(i_2)}\zeta_0^{(i_1)}-2\mu_q^{(i_1)}\zeta_0^{(i_2)}\Biggr)+
$$

\vspace{2mm}
$$
+\frac{1}{2\sqrt{2}}\sum_{r=1}^{q}
\Biggl(-\frac{1}{\pi r}
\zeta_{2r-1}^{(i_2)}
\zeta_{0}^{(i_1)}+
\frac{1}{\pi^2 r^2}\left(
\zeta_{2r}^{(i_2)}
\zeta_{0}^{(i_1)}-
2\zeta_{2r}^{(i_1)}
\zeta_{0}^{(i_2)}\right)\Biggr)-
$$

\vspace{2mm}
$$
-
\frac{1}{2\pi^2}\sum_{\stackrel{r,l=1}{{}_{r\ne l}}}^{q}
\frac{1}{r^2-l^2}\Biggl(
\zeta_{2r}^{(i_1)}
\zeta_{2l}^{(i_2)}+
\frac{l}{r}
\zeta_{2r-1}^{(i_1)}
\zeta_{2l-1}^{(i_2)}
\Biggr)+
$$

\vspace{2mm}
$$
+
\sum_{r=1}^{q}
\Biggl(\frac{1}{4\pi r}\left(
\zeta_{2r}^{(i_1)}
\zeta_{2r-1}^{(i_2)}-
\zeta_{2r-1}^{(i_1)}
\zeta_{2r}^{(i_2)}\right)+
$$

\vspace{2mm}
\begin{equation}
\label{944}
+
\Biggl.\Biggl.
\frac{1}{8\pi^2 r^2}\left(
3\zeta_{2r-1}^{(i_1)}
\zeta_{2r-1}^{(i_2)}+
\zeta_{2r}^{(i_2)}
\zeta_{2r}^{(i_1)}\right)\Biggr)\Biggr),
\end{equation}

\vspace{4mm}
\noindent
where the meaning of notations included 
in (\ref{42}), (\ref{46}) is saved.

An analogue of the formula (\ref{944}) (for the case of Legendre polynomials)
is (according to (\ref{4004}) and (\ref{4006})) the following representation

\vspace{2mm}
$$
I_{(10)T,t}^{*(i_1 i_2)q}=-\frac{T-t}{2}I_{(00)T,t}^{*(i_1 i_2)q}
-\frac{(T-t)^2}{4}\left(
\frac{\zeta_0^{(i_2)}\zeta_1^{(i_1)}}{\sqrt{3}}+\right.
$$

\vspace{3mm}
\begin{equation}
\label{5500}
+\left.\sum_{i=0}^{q}\left(
\frac{(i+1)\zeta_{i+2}^{(i_2)}\zeta_{i}^{(i_1)}
-(i+2)\zeta_{i}^{(i_2)}\zeta_{i+2}^{(i_1)}}
{\sqrt{(2i+1)(2i+5)}(2i+3)}+
\frac{\zeta_i^{(i_1)}\zeta_{i}^{(i_2)}}{(2i-1)(2i+3)}\right)\right),
\end{equation}

\vspace{7mm}
\noindent
where

\vspace{-1mm}
$$
I_{(00)T,t}^{*(i_1 i_2)q}=
\frac{T-t}{2}\left(\zeta_0^{(i_1)}\zeta_0^{(i_2)}+\sum_{i=1}^{q}
\frac{1}{\sqrt{4i^2-1}}\left(
\zeta_{i-1}^{(i_1)}\zeta_{i}^{(i_2)}-
\zeta_i^{(i_1)}\zeta_{i-1}^{(i_2)}\right)\right),
$$

\vspace{7mm}
\noindent
which is obviously substantially simpler than (\ref{944}).

Here it is necessary to pay a special attention on the fact
that the representation (\ref{5500}) includes a single sum with 
the upper summation limit $q$
while the representation (\ref{944}) includes the double sum 
with the same summation limit. In numerical simulation, 
obviously, the formula (\ref{5500}) is more economical 
in terms of computational cost than its analogue (\ref{944}).

There is another feature that should be noted
in connection with the formula (\ref{944}). This formula was first
obtained in \cite{KPW} by the method from \cite{Mi2}.
As we noted in Sect.~1, the method \cite{Mi2}
of approximation of iterated stochastic integrals 
is based on the trigonometric series expansion of the Brownian bridge 
process. So, this method 
leads to iterated application of the operation
of limit transition (in contrast to Theorems 1, 2, 6, and 8--12 in 
which limit transition is performed only once).
This means, generally speaking, that 
the mean-square convergence of 
$I_{(10)T,T}^{*(i_1 i_2)q}$
(see (\ref{944})) to 
$I_{(10)T,T}^{*(i_1 i_2)}$ does not follow if $q\to\infty$
for the method \cite{Mi2}.
The same applies to some others 
approximations of iterated Stratonovich stochastic integrals
obtained in \cite{KPW} by the method \cite{Mi2}
(see discussion in Sect.~8 for details).

The validity of the formula

$$
\lim\limits_{q\to\infty}
{\sf M}\left\{\left(I_{(10)T,t}^{*(i_1 i_2)}-I_{(10)T,t}^{*(i_1 i_2)q}
\right)^2\right\}=0,
$$

\vspace{5mm}
\noindent
where $I_{(10)T,T}^{*(i_1 i_2)q}$ is defined by (\ref{944}),
follows from Theorems 3, 6, and 8.

\vspace{5mm}

\section{Theorems 1, 2, 6, 8--12 from Point
of View of the Wong--Zakai Approximation}

\vspace{5mm}

The iterated Ito stochastic integrals and solutions
of Ito SDEs are complex and important func\-ti\-o\-nals
from the independent components ${\bf f}_{s}^{(i)},$
$i=1,\ldots,m$ of the multidimensional
Wiener process ${\bf f}_{s},$ $s\in[0, T].$
Let ${\bf f}_{s}^{(i)p},$ $p\in\mathbb{N}$ 
be some approximation of
${\bf f}_{s}^{(i)},$
$i=1,\ldots,m$.
Suppose that 
${\bf f}_{s}^{(i)p}$
converges to
${\bf f}_{s}^{(i)},$
$i=1,\ldots,m$ if $p\to\infty$ in some sense and has
differentiable sample trajectories.

A natural question arises: if we replace 
${\bf f}_{s}^{(i)}$
by ${\bf f}_{s}^{(i)p},$
$i=1,\ldots,m$ in the functionals
mentioned above, will the resulting
functionals converge to the original
functionals from the components 
${\bf f}_{s}^{(i)},$
$i=1,\ldots,m$ of the multidimentional
Wiener process ${\bf f}_{s}$?
The answere to this question is negative 
in the general case. However, 
in the pioneering works of Wong E. and Zakai M. \cite{W-Z-1},
\cite{W-Z-2},
it was shown that under the special conditions and 
for some types of approximations 
of the Wiener process the answere is affirmative
with one peculiarity: the convergence takes place 
to the iterated Stratonovich stochastic integrals
and solutions of Stratonovich SDEs and not to iterated 
Ito stochastic integrals and solutions
of Ito SDEs.
The piecewise 
linear approximation 
as well as the regularization by convolution 
\cite{W-Z-1}-\cite{Watanabe} relate the 
mentioned types of approximations
of the Wiener process. The above approximation 
of stochastic integrals and solutions of SDEs 
is often called the Wong--Zakai approximation.

Let ${\bf w}_{\tau},$ $\tau\in[0, T]$ is a random vector with 
an $m+1$ components: ${\bf w}_{\tau}^{(i)}={\bf f}_{\tau}^{(i)}$ 
for $i=1,\ldots,m$ and 
${\bf w}_{\tau}^{(0)}=\tau,$\ 
${\bf f}_{\tau}^{(i)}$ $(i=1,\ldots,m)$
are independent standard Wiener processes.

It is well known that the following representation 
takes place \cite{Lipt}, \cite{7e}

\vspace{-1mm}
\begin{equation}
\label{um1x}
{\bf w}_{\tau}^{(i)}-{\bf w}_{t}^{(i)}=
\sum_{j=0}^{\infty}\int\limits_t^{\tau}
\phi_j(s)ds\ \zeta_j^{(i)},\ \ \ \zeta_j^{(i)}=
\int\limits_t^T \phi_j(s)d{\bf w}_s^{(i)},
\end{equation}

\vspace{3mm}
\noindent
where $\tau\in[t, T],$ $t\ge 0,$
$\{\phi_j(x)\}_{j=0}^{\infty}$ is an arbitrary complete 
orthonormal system of functions in the space $L_2([t, T]),$ and
$\zeta_j^{(i)}$ are independent standard Gaussian 
random variables for various $i$ or $j.$
Moreover, the series (\ref{um1x}) converges for any $\tau\in [t, T]$
in the mean-square sense.

Let ${\bf w}_{\tau}^{(i)p}-{\bf w}_{t}^{(i)p}$ be 
the mean-square approximation of the process
${\bf w}_{\tau}^{(i)}-{\bf w}_{t}^{(i)},$
which has the following form

\vspace{-5mm}
\begin{equation}
\label{um1xx}
{\bf w}_{\tau}^{(i)p}-{\bf w}_{t}^{(i)p}=
\sum_{j=0}^{p}\int\limits_t^{\tau}
\phi_j(s)ds\ \zeta_j^{(i)}.
\end{equation}

\vspace{3mm}

From (\ref{um1xx}) we obtain

\vspace{-3mm}
\begin{equation}
\label{um1xxx}
d{\bf w}_{\tau}^{(i)p}=
\sum_{j=0}^{p}
\phi_j(\tau)\zeta_j^{(i)} d\tau.
\end{equation}

\vspace{4mm}

Consider the following iterated Riemann--Stieltjes
integral

\vspace{1mm}
\begin{equation}
\label{um1xxxx}
\int\limits_t^T
\psi_k(t_k)\ldots \int\limits_t^{t_2}\psi_1(t_1)
d{\bf w}_{t_1}^{(i_1)p_1}\ldots d{\bf w}_{t_k}^{(i_k)p_k},
\end{equation}

\vspace{4mm}
\noindent
where $p_1,\ldots,p_k\in\mathbb{N},$\ \ $i_1,\ldots,i_k=0,1,\ldots,m,$ 

\begin{equation}
\label{um1xxx1}
d{\bf w}_{\tau}^{(i)p}=
\left\{\begin{matrix}
d{\bf f}_{\tau}^{(i)p}\ &\hbox{\rm for}\ \ \ i=1,\ldots,m\cr\cr\cr
d\tau^p\ &\hbox{\rm for}\ \ \ i=0
\end{matrix}
,\right.
\end{equation}

\vspace{4mm}
\noindent
and $d{\bf f}_{\tau}^{(i)p},$ $d\tau^p$ are defined by the relation (\ref{um1xxx}).

Let us substitute (\ref{um1xxx}) into (\ref{um1xxxx})

\begin{equation}
\label{um1xxxx1}
\int\limits_t^T
\psi_k(t_k)\ldots \int\limits_t^{t_2}\psi_1(t_1)
d{\bf w}_{t_1}^{(i_1)p_1}\ldots d{\bf w}_{t_k}^{(i_k)p_k}=
\sum\limits_{j_1=0}^{p_1}\ldots \sum\limits_{j_k=0}^{p_k}
C_{j_k \ldots j_1}\prod\limits_{l=1}^k \zeta_{j_l}^{(i_l)},
\end{equation}

\vspace{4mm}
\noindent
where 
$$
\zeta_j^{(i)}=\int\limits_t^T \phi_j(s)d{\bf w}_s^{(i)}
$$ 

\vspace{2mm}
\noindent
are independent standard Gaussian random variables for various 
$i$ or $j$ (in the case when $i\ne 0$),
${\bf w}_{s}^{(i)}={\bf f}_{s}^{(i)}$ for
$i=1,\ldots,m$ and 
${\bf w}_{s}^{(0)}=s,$

$$
C_{j_k \ldots j_1}=\int\limits_t^T\psi_k(t_k)\phi_{j_k}(t_k)\ldots
\int\limits_t^{t_2}
\psi_1(t_1)\phi_{j_1}(t_1)
dt_1\ldots dt_k
$$

\vspace{4mm}
\noindent
is the Fourier coefficient.

To best of our knowledge \cite{W-Z-1}-\cite{Watanabe}
the approximations of the Wiener process
in the Wong--Zakai approximation must satisfy fairly strong
restrictions
\cite{Watanabe}
(see Definition 7.1, pp.~480--481).
Moreover, approximations of the Wiener process that are
similar to (\ref{um1xx})
were not considered in \cite{W-Z-1}, \cite{W-Z-2}
(also see \cite{Watanabe}, Theorems 7.1, 7.2).
Therefore, the proof of analogs of Theorems 7.1 and 7.2 \cite{Watanabe}
for approximations of the Wiener 
process based on its series expansion (\ref{um1x})
should be carried out separately.
Thus, the mean-square convergence of the right-hand side
of (\ref{um1xxxx1}) to the iterated Stratonovich stochastic integral 
(\ref{605})
does not follow from the results of the papers
\cite{W-Z-1}, \cite{W-Z-2} (also see \cite{Watanabe},
Theorems 7.1, 7.2).

From the other hand, Theorems 1, 2, 6, 8--12
can be considered as the proof of the
Wong--Zakai approximation for the iterated 
Stratonovich stochastic integrals (\ref{605}) of multiplicities 1--5 and $k$ ($k\in\mathbb{N}$)
based on the approximation (\ref{um1xx}) of the Wiener process.
At that, the Riemann--Stieltjes integrals (\ref{um1xxxx}) converge
(according to Theorems 1, 2, 6, 8--12)
to the appropriate Stratonovich 
stochastic integrals (\ref{605}). Recall that
$\{\phi_j(x)\}_{j=0}^{\infty}$ (see (\ref{um1x}), (\ref{um1xx}), and
Theorems 6, 9--12)
is a complete 
orthonormal system of Legendre polynomials or 
trigonometric functions 
in the space $L_2([t, T])$.

To illustrate the above reasoning, 
consider two examples for the case $k=2,$
$\psi_1(s),$ $\psi_2(s)\equiv 1;$ $i_1, i_2=1,\ldots,m.$

The first example relates to the piecewise linear approximation
of the multidimensional Wiener process (these approximations 
were considered in \cite{W-Z-1}-\cite{Watanabe}).

Let ${\bf b}_{\Delta}^{(i)}(t),$ $t\in[0, T]$ be the piecewise
linear approximation of the $i$th component ${\bf f}_t^{(i)}$
of the multidimensional standard Wiener process ${\bf f}_t,$
$t\in [0, T]$ with independent components
${\bf f}_t^{(i)},$ $i=1,\ldots,m,$ i.e.

\vspace{-2mm}
$$
{\bf b}_{\Delta}^{(i)}(t)={\bf f}_{k\Delta}^{(i)}+
\frac{t-k\Delta}{\Delta}\Delta{\bf f}_{k\Delta}^{(i)},
$$

\vspace{3mm}
\noindent
where 

\vspace{-1mm}
$$
\Delta{\bf f}_{k\Delta}^{(i)}={\bf f}_{(k+1)\Delta}^{(i)}-
{\bf f}_{k\Delta}^{(i)},\ \ \
t\in[k\Delta, (k+1)\Delta),\ \ \ k=0, 1,\ldots, N-1.
$$

\vspace{5mm}

Note that w.~p.~1

\vspace{-1mm}
\begin{equation}
\label{pridum}
\frac{d{\bf b}_{\Delta}^{(i)}}{dt}(t)=
\frac{\Delta{\bf f}_{k\Delta}^{(i)}}{\Delta},\ \ \
t\in[k\Delta, (k+1)\Delta),\ \ \ k=0, 1,\ldots, N-1.
\end{equation}

\vspace{4mm}

Consider the following iterated Riemann--Stieltjes
integral

\vspace{1mm}
$$
\int\limits_0^T
\int\limits_0^{s}
d{\bf b}_{\Delta}^{(i_1)}(\tau)d{\bf b}_{\Delta}^{(i_2)}(s),\ \ \ 
i_1,i_2=1,\ldots,m.
$$

\vspace{5mm}

Using (\ref{pridum})
and additive property of the Riemann--Stieltjes integral, 
we can write w.~p.~1

\vspace{2mm}
$$
\int\limits_0^T
\int\limits_0^{s}
d{\bf b}_{\Delta}^{(i_1)}(\tau)d{\bf b}_{\Delta}^{(i_2)}(s)=
\int\limits_0^T
\int\limits_0^{s}
\frac{d{\bf b}_{\Delta}^{(i_1)}}{d\tau}(\tau)d\tau
\frac{d {\bf b}_{\Delta}^{(i_2)}}{d s}(s)
ds =
$$

\vspace{3mm}
$$
=
\sum\limits_{l=0}^{N-1}\int\limits_{l\Delta}^{(l+1)\Delta}
\left(
\sum\limits_{q=0}^{l-1}\int\limits_{q\Delta}^{(q+1)\Delta}
\frac{\Delta{\bf f}_{q\Delta}^{(i_1)}}{\Delta}d\tau+
\int\limits_{l\Delta}^{s}
\frac{\Delta{\bf f}_{l\Delta}^{(i_1)}}{\Delta}d\tau\right)
\frac{\Delta{\bf f}_{l\Delta}^{(i_2)}}{\Delta}ds=
$$

\vspace{3mm}
$$
=\sum\limits_{l=0}^{N-1}\sum\limits_{q=0}^{l-1}
\Delta{\bf f}_{q\Delta}^{(i_1)}
\Delta{\bf f}_{l\Delta}^{(i_2)}+
\frac{1}{\Delta^2}\sum\limits_{l=0}^{N-1}
\Delta{\bf f}_{l\Delta}^{(i_1)}
\Delta{\bf f}_{l\Delta}^{(i_2)}
\int\limits_{l\Delta}^{(l+1)\Delta}
\int\limits_{l\Delta}^{s}d\tau ds=
$$

\vspace{3mm}
\begin{equation}
\label{oh-ty}
=\sum\limits_{l=0}^{N-1}\sum\limits_{q=0}^{l-1}
\Delta{\bf f}_{q\Delta}^{(i_1)}
\Delta{\bf f}_{l\Delta}^{(i_2)}+
\frac{1}{2}\sum\limits_{l=0}^{N-1}
\Delta{\bf f}_{l\Delta}^{(i_1)}
\Delta{\bf f}_{l\Delta}^{(i_2)}.
\end{equation}

\vspace{6mm}

Using (\ref{oh-ty}) and  (\ref{oop51}) it 
is not difficult to show 
that

\vspace{1mm}
$$
\hbox{\vtop{\offinterlineskip\halign{
\hfil#\hfil\cr
{\rm l.i.m.}\cr
$\stackrel{}{{}_{N\to \infty}}$\cr
}} }
\int\limits_0^T
\int\limits_0^{s}
d{\bf b}_{\Delta}^{(i_1)}(\tau)d{\bf b}_{\Delta}^{(i_2)}(s)=
\int\limits_0^T
\int\limits_0^{s}
d{\bf f}_{\tau}^{(i_1)}d{\bf f}_{s}^{(i_2)}+
\frac{1}{2}{\bf 1}_{\{i_1=i_2\}}\int\limits_0^T ds=
$$

\vspace{3mm}
\begin{equation}
\label{uh-111}
=
{\int\limits_0^{*}}^{T}
{\int\limits_0^{*}}^{s}
d{\bf f}_{\tau}^{(i_1)}d{\bf f}_{s}^{(i_2)},
\end{equation}

\vspace{5mm}
\noindent
where $\Delta\to 0$ if $N\to\infty$ ($N\Delta=T$).

Obviously, (\ref{uh-111}) agrees with Theorem 7.1 (see \cite{Watanabe},
p.~486).

The next example relates to the approximation
of the Wiener process based on its series expansion
(\ref{um1x}) for $t=0$, where
$\{\phi_j(x)\}_{j=0}^{\infty}$ 
is a complete 
orthonormal system of Legendre polynomials or 
trigonometric functions 
in the space $L_2([0, T])$.

Consider the following iterated Riemann--Stieltjes
integral

\vspace{-1mm}
\begin{equation}
\label{abcd1}
\int\limits_0^T
\int\limits_0^{s}
d{\bf f}_{\tau}^{(i_1)p}d{\bf f}_{s}^{(i_2)p},\ \ \ 
i_1,i_2=1,\ldots,m,
\end{equation}

\vspace{3mm}
\noindent
where $d{\bf f}_{\tau}^{(i)p}$ is defined by the
relation
(\ref{um1xxx}).

Let us substitute (\ref{um1xxx}) into (\ref{abcd1}) 

\vspace{-1mm}
\begin{equation}
\label{set18}
\int\limits_0^T
\int\limits_0^{s}
d{\bf f}_{\tau}^{(i_1)p}d{\bf f}_{s}^{(i_2)p}=
\sum\limits_{j_1,j_2=0}^p
C_{j_2 j_1} \zeta_{j_1}^{(i_1)}\zeta_{j_2}^{(i_2)},
\end{equation}

\vspace{3mm}
\noindent
where 
$$
C_{j_2 j_1}=
\int\limits_0^T \phi_{j_2}(s)\int\limits_0^s
\phi_{j_1}(\tau)d\tau ds
$$

\vspace{3mm}
\noindent
is the Fourier coefficient; another notations 
are the same as in (\ref{um1xxxx1}).

As we noted above, approximations of the Wiener process that are
similar to (\ref{um1xx})
were not considered in \cite{W-Z-1}, \cite{W-Z-2}
(also see Theorems 7.1, 7.2 in \cite{Watanabe}).
Furthermore, the extension of the results of Theorems 7.1 and 7.2
\cite{Watanabe} to the case under consideration is
not obvious.

However, the authors of the works
\cite{KlPl2}
(Sect.~5.8, pp.~202--204), \cite{KPW} (pp.~438-439),  
\cite{KPS} (pp.~82-84),
\cite{Zapad-9} (pp.~263-264) use 
the Wong--Zakai approximation 
\cite{W-Z-1}-\cite{Watanabe} (without rigorous proof) within the frames
of the approach
\cite{Mi2} based on the series expansion 
of the Brownian bridge process.

On the other hand, we can apply the theory built in Chapters 1 and 2
of the monographs \cite{20xxxxx}-\cite{20eee}. More precisely, 
using 
Theorem 6 from this paper 
we obtain from (\ref{set18}) the desired result

\vspace{-1mm}
$$
\hbox{\vtop{\offinterlineskip\halign{
\hfil#\hfil\cr
{\rm l.i.m.}\cr
$\stackrel{}{{}_{p\to \infty}}$\cr
}} }
\int\limits_0^T
\int\limits_0^{s}
d{\bf f}_{\tau}^{(i_1)p}d{\bf f}_{s}^{(i_2)p}=
\hbox{\vtop{\offinterlineskip\halign{
\hfil#\hfil\cr
{\rm l.i.m.}\cr
$\stackrel{}{{}_{p\to \infty}}$\cr
}} }
\sum\limits_{j_1,j_2=0}^p
C_{j_2 j_1} \zeta_{j_1}^{(i_1)}\zeta_{j_2}^{(i_2)}=
$$

\vspace{2mm}
\begin{equation}
\label{umen-bl}
=
{\int\limits_0^{*}}^{T}
{\int\limits_0^{*}}^{s}
d{\bf f}_{\tau}^{(i_1)}d{\bf f}_{s}^{(i_2)}.
\end{equation}

\vspace{4mm}

From the other hand, by Theorems 1, 2
(see (\ref{leto5001})) for the case
$k=2$ we obtain from (\ref{set18}) the following relation

\vspace{-2mm}
$$
\hbox{\vtop{\offinterlineskip\halign{
\hfil#\hfil\cr
{\rm l.i.m.}\cr
$\stackrel{}{{}_{p\to \infty}}$\cr
}} }
\int\limits_0^T
\int\limits_0^{s}
d{\bf f}_{\tau}^{(i_1)p}d{\bf f}_{s}^{(i_2)p}=
\hbox{\vtop{\offinterlineskip\halign{
\hfil#\hfil\cr
{\rm l.i.m.}\cr
$\stackrel{}{{}_{p\to \infty}}$\cr
}} }
\sum\limits_{j_1,j_2=0}^p
C_{j_2 j_1} \zeta_{j_1}^{(i_1)}\zeta_{j_2}^{(i_2)}=
$$

\vspace{2mm}
$$
=
\hbox{\vtop{\offinterlineskip\halign{
\hfil#\hfil\cr
{\rm l.i.m.}\cr
$\stackrel{}{{}_{p\to \infty}}$\cr
}} }
\sum\limits_{j_1,j_2=0}^p
C_{j_2 j_1} \biggl(\zeta_{j_1}^{(i_1)}\zeta_{j_2}^{(i_2)}-
{\bf 1}_{\{i_1=i_2\}}{\bf 1}_{\{j_1=j_2\}}\biggr)+
{\bf 1}_{\{i_1=i_2\}}\sum\limits_{j_1=0}^{\infty}
C_{j_1 j_1}=
$$

\vspace{2mm}
\begin{equation}
\label{umen-blx}
=
\int\limits_0^T
\int\limits_0^{s}
d{\bf f}_{\tau}^{(i_1)}d{\bf f}_{s}^{(i_2)}+
{\bf 1}_{\{i_1=i_2\}}\sum\limits_{j_1=0}^{\infty}
C_{j_1 j_1}.
\end{equation}

\vspace{5mm}

Since
$$
\sum\limits_{j_1=0}^{\infty}
C_{j_1 j_1}=\frac{1}{2}\sum\limits_{j_1=0}^{\infty}
\left(\int\limits_0^T \phi_j(\tau)d\tau\right)^2
=
$$

\vspace{2mm}
$$
=\frac{1}{2}
\left(\int\limits_0^T \phi_0(\tau)d\tau\right)^2
=\frac{1}{2}
\int\limits_0^T ds,
$$

\vspace{5mm}
\noindent
then from (\ref{oop51}) and (\ref{umen-blx}) we obtain (\ref{umen-bl}).

\end{document}